\documentclass[12pt, twoside]{article}
\usepackage{amsmath,amsthm,amssymb}
\usepackage{times}
\usepackage{enumerate}

\usepackage{amsfonts,amscd,epsfig,psfrag}
\usepackage{graphicx}
\usepackage{yhmath}
\usepackage[mathscr]{eucal}
\usepackage{xspace}
\newcommand{\R}{\mathbb{R}}
\newcommand{\D}{\mathbb{D}}

\newcommand{\N}{\mathbb{N}}

\newcommand{\CS}{\mathcal{S}}

\newcommand{\xx}{\mathbf{x}}
\newcommand{\zz}{\mathbf{z}}
\newcommand{\ff}{\mathbf{f}}
\newcommand{\kk}{\mathbf{k}}

\newcommand{\uu}{\mathbf{u}}

\newcommand{\ga}{\alpha}
\newcommand{\gb}{\mathbf{\beta}}
\newcommand{\gc}{\gamma}

\newcommand{\gD}{\Delta}
\newcommand{\po}{\partial}

\newcommand{\ve}{\varepsilon}
\newcommand{\vae}{\varepsilon}

\newcommand{\gd}{\delta}

\newcommand{\gl}{\lambda}

\newcommand{\del}{{\partial}}
\renewcommand{\O}{{\mathbf O}}
\renewcommand{\d}{\delta}
\renewcommand{\l}{\lambda}

\renewcommand{\b}{\beta}
\renewcommand{\t}{\theta}
\newcommand{\s}{\sigma}
\newcommand{\g}{\gamma}

\newcommand{\Om}{\Omega}
\newcommand{\om}{\omega}

\newcommand{\grad}{\nabla}

\newcommand{\tU}{\widetilde{U}}

\newcommand{\tv}{\tilde{v}}
\newcommand{\tp}{\tilde{p}}

\newcommand{\bea}{\begin{eqnarray}}
\newcommand{\eea}{\end{eqnarray}}

\newcommand{\yy}{\mathbf{y}}

\newcommand{\hs}{\hat{\sigma}}
\newcommand{\pw}{\partial \mathcal{W}}
\newcommand{\tgs}{\widetilde{\Sigma}}

\theoremstyle{plain}
\newtheorem{theorem}{Theorem}[section]

\newtheorem*{problem}{Problem}
\newtheorem{lemma}{Lemma}[section]

\theoremstyle{definition}

\theoremstyle{remark}

\newtheorem{remark}{Remark}[section]

\numberwithin{equation}{section}
\numberwithin{figure}{section}
\newcommand{\subjclass}[1]{\bigskip\noindent\emph{2010 Mathematics Subject Classification:}\enspace#1}
\newcommand{\keywords}[1]{\noindent\emph{Keywords:}\enspace#1}
\begin{document}

%%%%% To ease editing, add:

\baselineskip=17pt

%%%%%%%%%%%%%%%%

\title{Stability and Asymptotic Behavior
of Transonic Flows Past Wedges for the Full Euler Equations}

\author{Gui-Qiang G. Chen\\
{\em Mathematical Institute, University of Oxford, Oxford,  OX2 6GG, UK}\\
{\em AMSS \& UCAS, Chinese Academy of Sciences, Beijing 100190, China}\\
{\em \small E-mail:  chengq@maths.ox.ac.uk}\\
Jun Chen\\
{\em Department of Mathematics, Southern University of Science and Technology}\\
{\em Shenzhen, Guangdong 518055, China}\\
{\em \small E-mail: chenjun@sustc.edu.cn}\\
Mikhail Feldman\\
{\em Department of Mathematics, University of Wisconsin-Madison}\\
{\em Madison, WI 53706-1388, USA}\\
{\em \small E-mail: feldman@math.wisc.edu}}

\date{[Received 5 February 2017]}

\maketitle

\begin{abstract}
The existence, uniqueness, and asymptotic behavior of
steady transonic flows past a curved wedge, involving transonic shocks,
governed by the two-dimensional full Euler equations  are established.
The stability of both weak and strong transonic shocks under the perturbation
of the upstream supersonic flow and the wedge boundary is proved.
The problem is formulated as a one-phase free boundary problem,
in which the transonic shock is treated as a free boundary.
The full Euler equations are decomposed into two algebraic equations
and a first-order elliptic system of two equations in Lagrangian coordinates.
With careful elliptic estimates by using appropriate weighted H\"older norms,
the iteration map is defined and analyzed, and the existence of its fixed point is established
by performing the Schauder fixed point argument.
The careful analysis of the asymptotic behavior of the solutions reveals
particular characters of the full Euler equations.

\subjclass{Primary: 35R35, 35M12,
  76H05, 76L05,
  35L67, 35L65,
  35B35,  35B30, 35B40, 35Q31,
  76N10, 76N15,
  35L60;
Secondary: 35M10,
  35B65,
  35B45, 35J67,
   76J20, 76G25.}

\keywords{Shock wave, free boundary, wedge problem, steady,
 supersonic, subsonic, transonic, mixed type, composite type,
 hyperbolic-elliptic, full Euler equations, physical admissible, existence, stability,
 asymptotic behavior, decay rate.}
\end{abstract}

\section{Introduction}

\begin{figure}
 	\centering
 	\includegraphics[height=50mm]{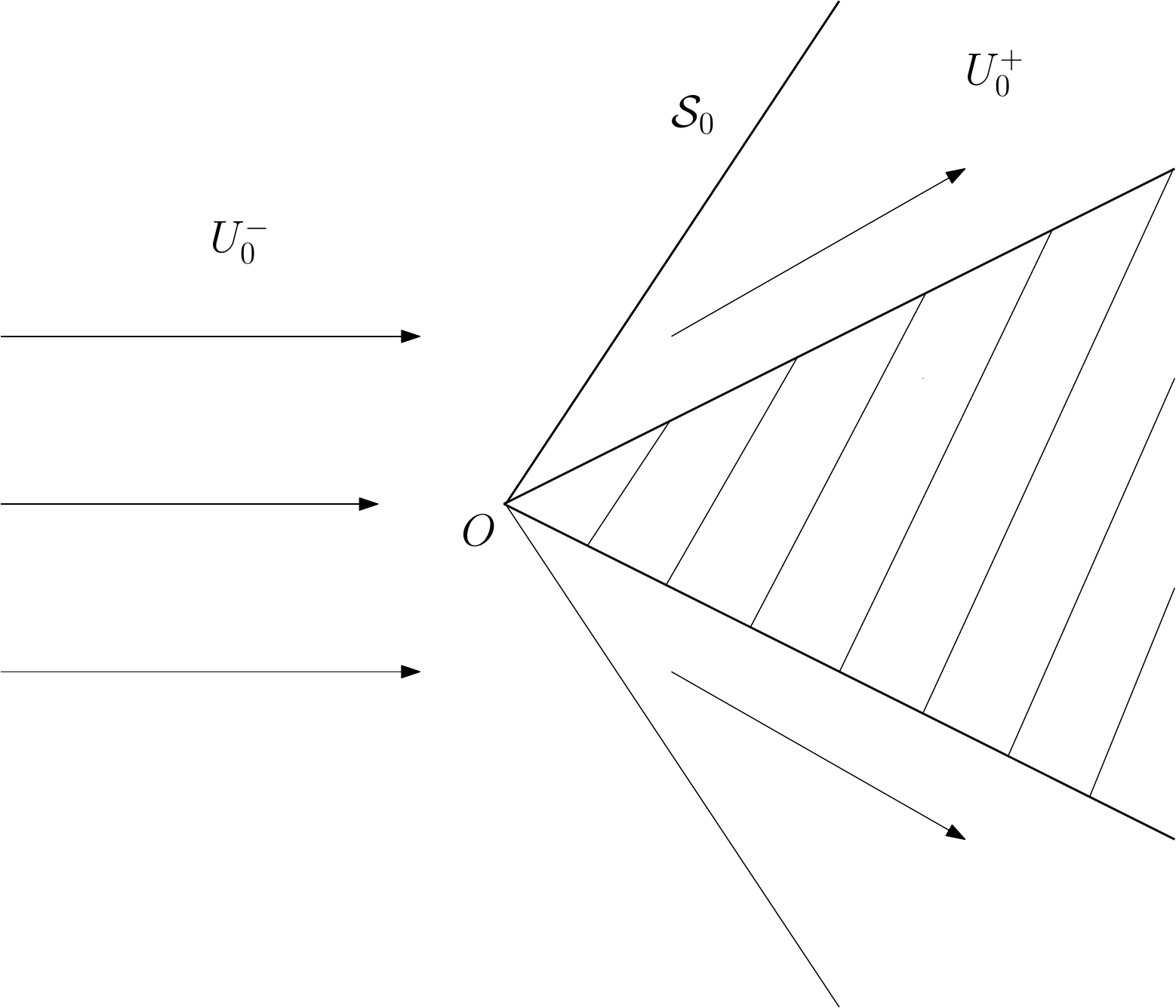}
\hspace{15mm}
 		\includegraphics[height=45mm]{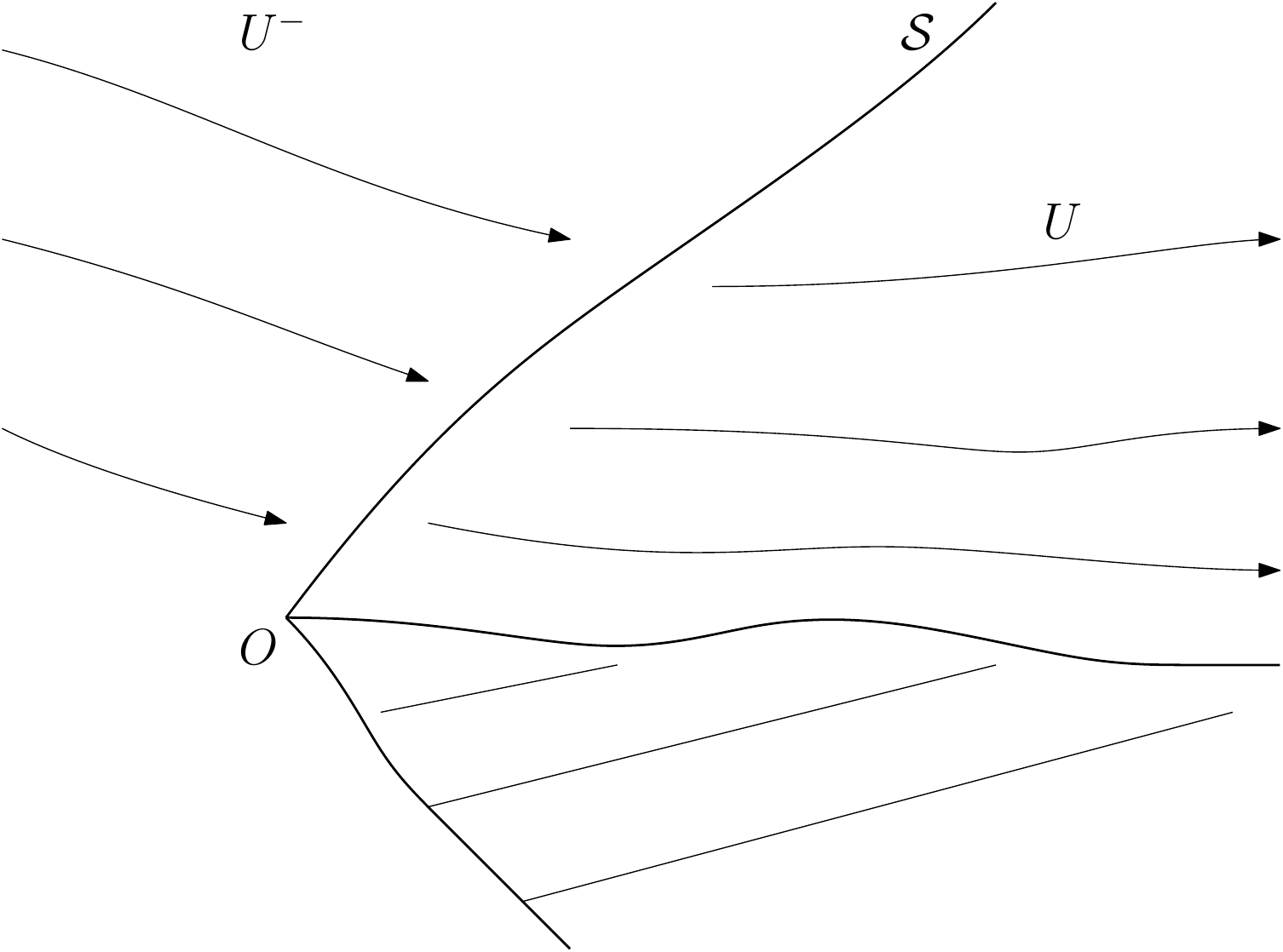}
		\caption{\small (a) Constant transonic flows; (b) Perturbed transonic flows rotated clockwise by angle $\theta_{\rm w}$}
		\label{fig-perturbedshock}
 \end{figure}

\begin{figure}
 \centering
\includegraphics[height=68mm]{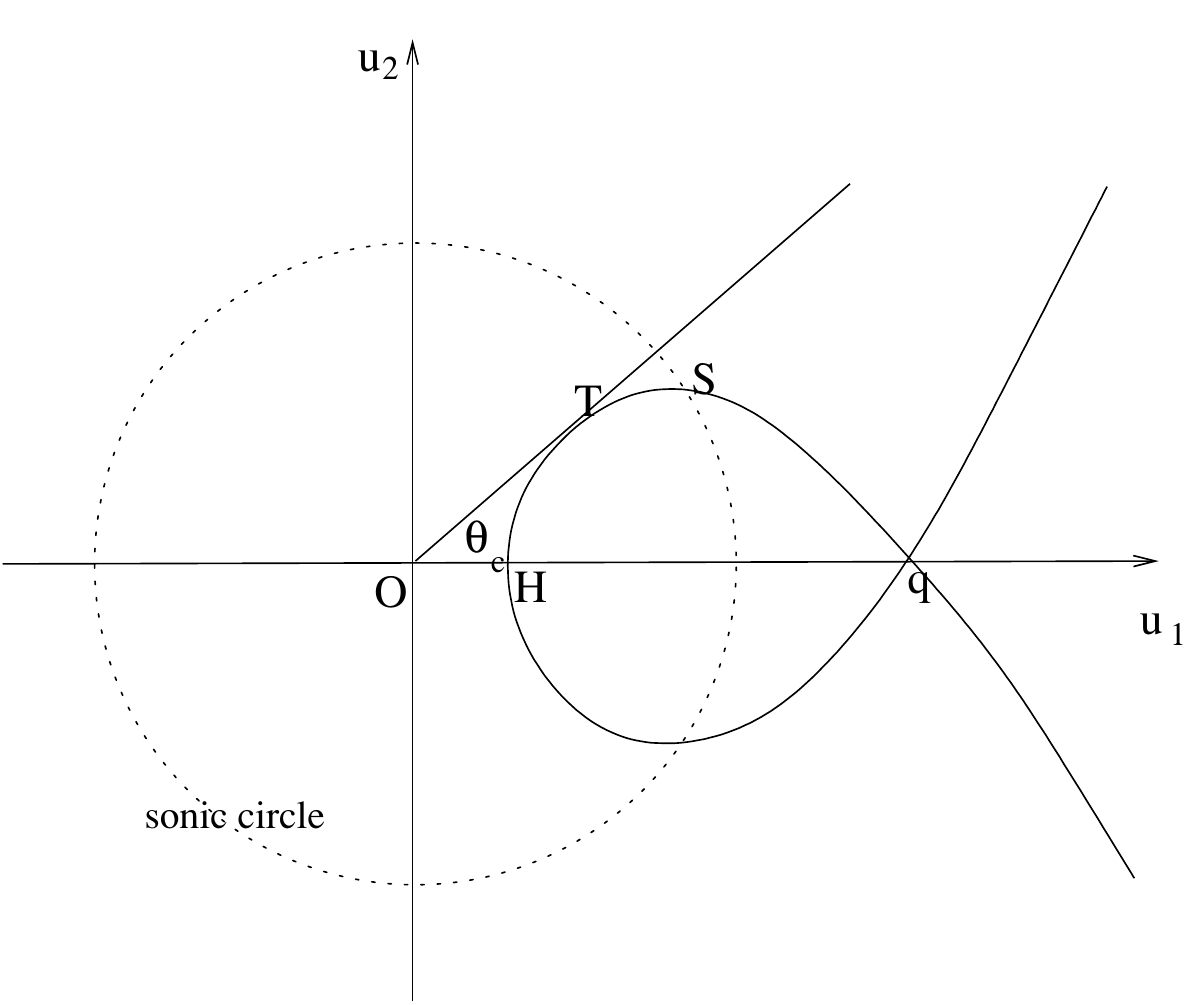}
\caption{The shock polar in the $(u_1, u_2)$--plane}
\label{fig-shockpolaruv}
\end{figure}

We are concerned with the existence, uniqueness, and asymptotic behavior of
steady transonic flows past a curved wedge, involving transonic shocks,
governed by the two-dimensional full Euler equations.
When a supersonic flow passes through a straight-sided wedge whose half-angle $\theta_{\rm w}$ is less than
the detachment angle, a shock attached to the wedge vertex is expected to form.
If the upstream steady flow is a uniform supersonic state, we can
find the corresponding
constant flow downstream along the straight-sided wedge boundary,
together with a straight shock separating
the two states (see Fig.~\ref{fig-perturbedshock}),
by using the shock polar determined by the Rankine-Hugoniot jump conditions
and the entropy condition ({\it cf.} Fig.~\ref{fig-shockpolaruv}).
However, these conditions do not determine the downstream state uniquely.
In general, there are two solutions, one of which corresponds to a weaker shock than the other.
As normally expected, a physically admissible shock should be stable under small perturbations.
Therefore, it is important to analyze the stability of these shocks
in order to understand underlying physics.

The wedge problem described above has a long history at least dating back to the 1930s.
Prandtl \cite{Prandtl} in 1936 first conjectured that the weak shock solution is stable,
and hence
physically admissible.
There has been a long debate about whether the strong shock is stable
for decades; see Courant-Friedrichs \cite{CourantF}, Section 123,
and von Neumann \cite{Neumann}. See also Liu \cite{Liu} and Serre \cite{Serre}.

When the downstream flow is supersonic, the corresponding shock is
called a supersonic shock,
which is a weak shock. This case has been analyzed for the potential flow
equation in \cite{Chen1,Chen2} with certain convexity assumption on the wedge
and in \cite{Zh2} for an almost straight-sided wedge.
The existence and stability of the steady supersonic shocks for the full Euler equations
have been established under the BV perturbation of both the upstream flow
and the slope of the wedge boundary
in Chen-Zhang-Zhu \cite{ChenZhangZhu} and Chen-Li \cite{ChenLi} for
Lipschitz wedges.

For transonic
shocks ({\it i.e.}, the downstream flow is subsonic),
there are two cases: the transonic shock with the subsonic state corresponding
to arc $\wideparen{TS}$ (which is a weak shock) and the one corresponding
to arc $\wideparen{TH}$ (which is a strong shock) (see Fig.~\ref{fig-shockpolaruv}).
The strong shock case has been studied in Chen-Fang \cite{ChenFang} for the potential flow
(also see \cite{CFang2}).

It is well known that the jump of the entropy function across the shock is of cubic order of
the shock strength.
In general, the strength of transonic shocks is large, so the full Euler system
is a more accurate model than the potential flow or isentropic Euler equations.
In Fang \cite{Fang},  the Euler equations were first studied with a uniform Bernoulli constant
for both strong and weak transonic shocks. However, the asymptotic behavior of the shock slope or the subsonic
part of the solution was not analyzed in \cite{Fang},
partly because the approach in \cite{Fang} is based on
the weighted Sobolev spaces.
On the other hand, the asymptotic behavior can be seen more conveniently within the framework of  H\"older spaces.
In Yin-Zhou \cite{YinZhou}, the H\"older norms were used for the estimates of the full Euler equations
with the assumption on the sharpness of the wedge angle, which means that the subsonic state
is near point $H$ in the shock polar.
In Chen-Chen-Feldman \cite{CCF2}, the weak transonic shock, which corresponds to
the whole arc $\wideparen{TS}$, was investigated; and
the existence, uniqueness, stability, and asymptotic
behavior of  subsonic solutions were obtained.
In \cite{CCF2,YinZhou}, a potential function is used to reduce the four Euler equations into
one elliptic equation in the subsonic region. The method was first proposed in \cite{CCF} and
has the advantage of integrating the conservation properties of the Euler system
into a single elliptic
equation. However, working on the potential function further requires its Lipschitz estimate, besides
the $C^0$--estimate, to keep the subsonicity of the flow.

There are other related papers about transonic shocks,
such as \cite{CFang2,LiangXuYin} for transonic flows past three-dimensional wedges
and \cite{CFang1} about transonic flows  past a perturbed cone;
see also \cite{CFeldman1,schen3} for the approaches developed earlier for dealing with
transonic shock flows and \cite{FangXiang} for the uniqueness of transonic shocks.

The purpose of this paper is to analyze both strong and weak transonic shocks and
establish the existence, uniqueness, and asymptotic behavior of the subsonic solutions
under the perturbation of both the upstream supersonic flows and the wedge boundaries.
In particular, we are able to prove the stability of both weak and strong
transonic shocks.
The strategy is to use the physical variables to make the estimates, instead of the potential function.
The advantage of this method is that only the lower regularity ({\it i.e.}, the $C^0$--estimate)
is sufficient to guarantee the subsonicity.
Furthermore, estimating the physical state function $U=(\uu, p, \rho)^\top$
directly (see equations \eqref{Euler1})
also yields a better asymptotic decay rate:
For weak transonic shocks,  the decay rate is only $|\xx|^{-\gb}$ in our earlier paper \cite{CCF2};
while, in this paper, we will show that the subsonic solution decays to a limit state at
rate, $|\xx|^{-1-\gb}$, with $\beta\in (0,1)$ depending only on the background states (see Remark \ref{rem:2.2a}).

More precisely, we first use the Lagrangian coordinates to straighten the streamlines.
The reason for this is that the Bernoulli variable and entropy are conserved along the streamlines,
and using the streamline as one of the coordinates simplifies the formulation,
especially for the asymptotic behavior of the solution.
Then, as in \cite{schen3,Fang},
we decompose the Euler system into two algebraic equations and two elliptic equations.
Differentiating the two elliptic equations gives rise to a second-order elliptic equation
in divergence form for the flow direction $w= \frac{u_2}{u_1}$.
Given $U$ in an expected function space for solutions, we obtain the updated function
$\tilde{w}$
as the solution of the linear equation for iterations
whose coefficients are evaluated on the given function $U$.
Once we solve for $\tilde{w}$ and obtain the desired estimates,
the other variables are then updated.
Thus, we construct a map $\d \widetilde{U}=\mathcal{Q}(\d U)$,
where $\d U$ and $\d \widetilde{U}$ are the perturbations from the background subsonic state.
The estimates based on our method do not yield the contraction for $\mathcal{Q}$.
Therefore, the Banach fixed point argument does not work; Instead, we employ the Schauder
fixed point argument to obtain the existence of the subsonic
solution. For the uniqueness, we estimate the difference of two solutions
by using the weighted H\"older norms with a lower decay rate.

One point we want to emphasize here is that the decay pattern is different from that for potential flow.
In a potential flow, the decay is with respect to $|\xx|$.
For example, if $\varphi$ converges to $\varphi_0$ at rate $|\xx|^{-\gb}$,
then $\grad \varphi$ converges at rate $|\xx|^{-1-\gb}$.
For the Euler equations, because the Bernoulli variable and the entropy function
are constant along streamlines,
the physical variables $(u_1, \rho)$ do not converge to the background state along the streamlines.
They converge only across the streamlines away from the wedge.
Therefore, when the elliptic estimates are performed,
the scaling is with respect to the distance from the wedge,
rather than  $|\xx|$.
This results in the following decay pattern:
In Lagrangian coordinates $\yy$, there exists an asymptotic limit $U^\infty = (u_1^\infty, 0, p_0^+, \rho^\infty)$;
$U$ converges to $U^\infty$ at rate $|\yy|^{-\gb}$, but $\grad U$ converges at rate $|\yy|^{-\gb} (y_2+1)^{-1}$.
That is, the extra decay for the derivatives is only along the $y_2$--direction.

Finally, we remark that our analysis of transonic shocks for the Euler equations
for potential and non-potential flows,
started in Chen-Feldman \cite{CFeldman1} to formulate the transonic shock problems
as  one-phase free boundary problems,
is motivated by the previous works on variational one-phase free boundary problems
for nonlinear elliptic equations in Alt-Caffarelli \cite{ACa},
Alt-Caffarelli-Friedman \cite{ACF1,ACF2},
and the references cited therein.
One of the main difficulties in dealing with the transonic shock problems is that the corresponding
elliptic one-phase free boundary problems are non-variational in general,
so that the complete solution to
the free boundary problems
requires different approaches and new techniques which are further developed in this paper
in the physical realm of the full Euler equations for compressible fluids.

The rest of the paper is organized in the following sections.
In \S \ref{sec-setup}, the wedge problem is formulated as a free boundary problem and
the main theorem is stated.
In \S \ref{sec-lagrange}, the problem is reformulated in Lagrangian coordinates.
In \S \ref{sec-decompose},
the Euler equations are decomposed into two algebraic equations and a first-order elliptic system of two equations.
In \S \ref{sec-linearize}, the linear elliptic system and the boundary conditions for iterations are introduced.
In \S \ref{sec-keylemma}, the key estimates of solutions for the linear second-order elliptic equation for iterations are obtained.
In \S \ref{sec-iteration}, the iteration map is constructed and the corresponding estimates are obtained,
leading to the existence of a weak transonic shock solution.
In \S \ref{sec-unique}, the uniqueness of the weak transonic shock solution is  proved.
In \S \ref{sec-asymptotic}, the asymptotic behavior and the decay rate of solutions are discussed.
In \S \ref{sec-th}, the difference between the weak and the strong transonic shocks
is revealed in terms of the estimates and the asymptotic behavior of the solution.

\section{Mathematical Setup and the Main Theorem} \label{sec-setup}

In this section, we formulate the transonic wedge problem as a
free boundary problem  and state the main theorem.

The governing equations are two-dimensional steady, full
 Euler equations:
 \begin{equation}
 \left\{\begin{aligned}
 & \nabla\cdot (\rho \uu)=0, \\
 &\nabla\cdot\left(\rho{\uu\otimes\uu}\right)
 +\nabla p=0,\\
 &\nabla\cdot\big(\rho\uu(E+\frac{p}{\rho})\big)=0,
 \end{aligned}\right.
 \label{Euler1}
 \end{equation}
 where $\grad$ is the gradient in $\xx=(x_1,x_2)\in\R^2$,
 $\uu=(u_1,  u_2)$ the velocity, $\rho$ the density, $p$ the pressure,
 and $\gamma>1$ the adiabatic exponent,
 as well as
 $$
 E=\frac{1}{2}|\uu|^2+\frac{p}{(\gamma-1)\rho}
 $$
 is the energy.
 The sonic speed of  the flow is
 $$
 c=\sqrt{\frac{\gamma p}{\rho}}.
 $$
 The flow is subsonic if $|\uu| < c$
 and supersonic if $|\uu| > c$.
 For a transonic flow, both cases occur in the flow.

System \eqref{Euler1} can be written in the following general form as
a system of conservation laws:
\begin{equation}\label{Euler1a}
\nabla\cdot \mathbf{F}(U)=0,  \qquad\,\, \xx\in \R^2,
\end{equation}
with $U=(\uu, p,\rho)^\top$.
Such systems
often govern time-independent solutions for multidimensional quasilinear
hyperbolic systems of conservation laws; {\it cf.} Dafermos \cite{Dafermos}
and Lax \cite{Lax}.

To be a weak
solution of the Euler equations \eqref{Euler1}, the Rankine-Huguoniot
conditions must be satisfied along the shock-front $x_1=\sigma(x_2)$:
  \begin{equation} \label{con-RH}
  \left\{ \begin{aligned}
  &[\,\rho u_1\,]=\s'(x_2) [\,\rho u_2\,],\\
  &[\,\rho u_1^2 + p\,]= \s'(x_2) [\,\rho u_1 u_2\,],\\
  &[\,\rho u_1 u_2\,] = \s'(x_2) [\,\rho {u_2}^2 + p\,],\\
  &[\,\rho u_1(E+\frac{p}{\rho})\,] = \s'(x_2) [\,\rho u_2 (E+\frac{p}{\rho})\,],
  \end{aligned} \right.
  \end{equation}
 where $[\,\cdot\,]$ denotes the jump of the quantity between the two states
 across the shock front; that is,
 if  $w^-$ and $w^+$ represent the left and right states, respectively,
 then $[w]:= w^+ - w^-$.

For a given constant upstream supersonic flow
$U_0^-=(u_{10}^-,0, p_0^-, \rho_0^-)^\top$
and a fixed straight-sided wedge with wedge angle $\theta_{\rm w}$,
the downstream constant flow
can be determined by the Rankine-Huguoniot conditions \eqref{con-RH1}--\eqref{con-RH4}.
According to the shock polar (see Fig.~\ref{fig-shockpolaruv}),
there are two subsonic solutions (for a large-angle wedge),
or one subsonic solution and one supersonic solution (for a small-angle wedge).
We choose the subsonic constant state for the downstream flows.
When the wedge angle $\theta_{\rm w}$ is between $0$ and the detachment
angle $\theta_{\rm w}^{\rm d}$,
arc $ \wideparen{HS}$ is divided by the tangent point $T$
into two open arcs $\wideparen{TH}$ and $\wideparen{TS}$,
which correspond to the strong and weak transonic shocks,
respectively.

For convenience, we rotate the plane clockwise by angle $\theta_{\rm w}$ so that the downstream
flows become horizontal.
Then $\frac{u_{20}^-}{u_{10}^-}=-\tan\theta_{\rm w}$,
$U_0^-=(u_{10}^-,-u_{10}^-\tan \theta_{\rm w}, p_0^-, \rho_0^-)^\top$,
and $U_0^+ =(u_{10}^+,0, p_0^+, \rho_0^+)^\top$ ({\it cf.} Fig.~\ref{fig-perturbedshock}).

Suppose that the background shock is the straight line
given by $\mathcal{S}_0:=\{x_1= \s_0(x_2):=k_0 x_2\}$.
Let $\Om^- $ be the region for the upstream flows defined by
\begin{equation*}
\Om^- = \Big\{\xx\,:\, 0<x_1< \frac{4}{3}k_0 x_2\Big\}.
\end{equation*}

We use a function $b(x_1)$ to describe the wedge boundary:
\begin{equation}\label{wall}
\del \mathcal{W}:=\{\xx\in \R^2\,:\,  x_2= b(x_1),\ b(0)=0\}.
\end{equation}
Along the solid wedge boundary $\del\mathcal{W}$, the slip condition is satisfied:
\begin{equation}\label{slipcon}
\left. \frac{u_2}{u_1}\right|_{\del \mathcal{W}} = b'.
\end{equation}

Suppose that the shock front $\mathcal{S}$ we seek is
\[
\mathcal{S}:=\{\xx \,:\, \s(0)=0, \, x_1 = \s(x_2), x_2 \ge 0\}.
\]
Then the domain for the subsonic flow is denoted by
\begin{equation}\label{domain:1}
\Omega^\s:=\{\xx\in \R^2\,:\,  x_1> \s(x_2),\,  x_2 > b(x_1)\}.
\end{equation}

\begin{figure}
 \centering
\includegraphics[height=58mm]{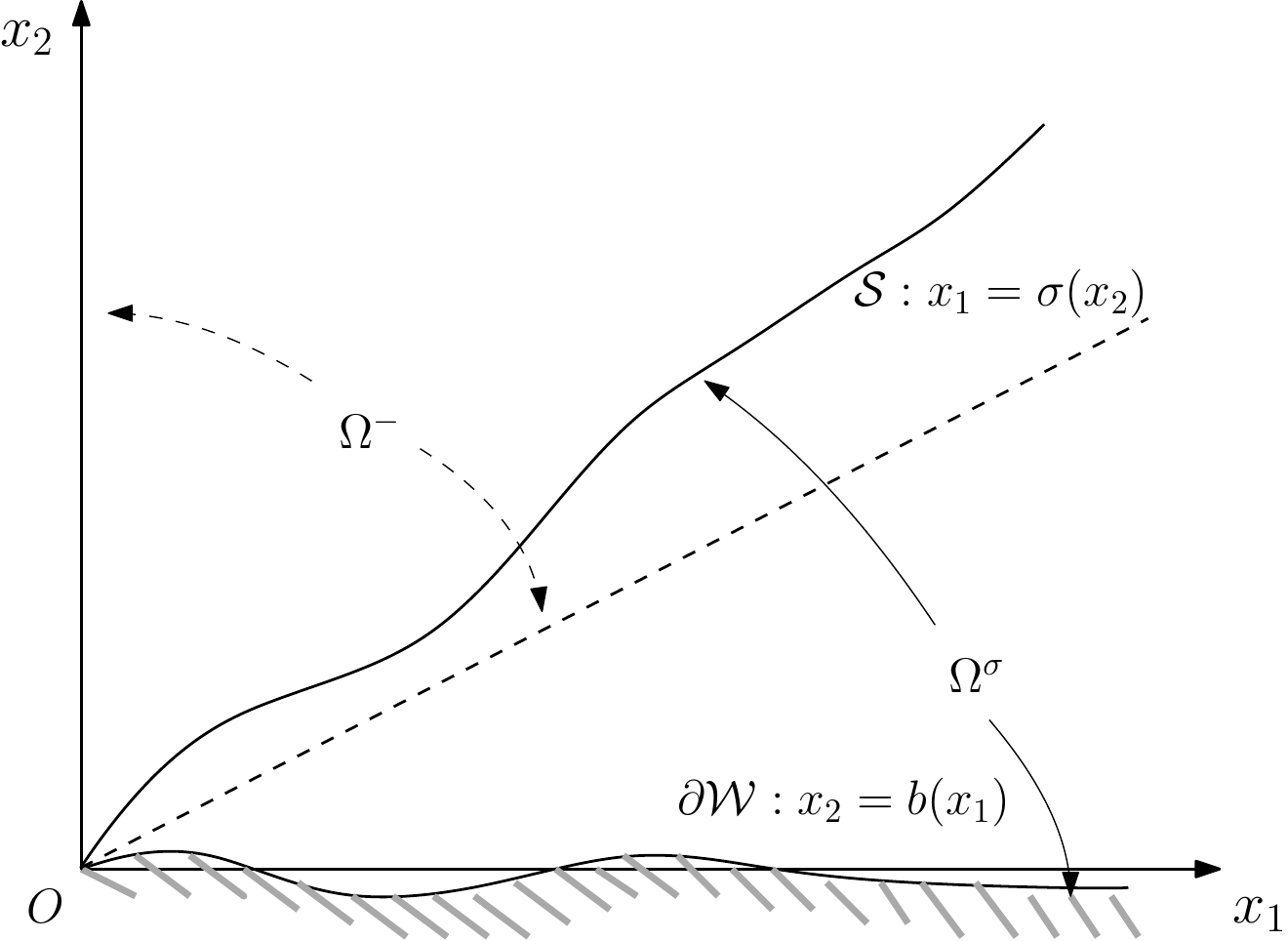}
\caption{Domains $\Omega^{-}$ and $\Omega^\s$ in Eulerian coordinates}
\label{new-1}
\end{figure}

Therefore, the problem can be formulated as the following free boundary problem:

\begin{problem}[Free Boundary Problem; see Fig. \ref{new-1}]
Let $(U^-_0, U^+_0)$ be a constant transonic solution with transonic shock $\mathcal{S}_0$.
For any upstream flow $U^-$ for equations \eqref{Euler1} in domain $\Om^-$,
which is a small perturbation of $U^-_0$,
find a subsonic solution $U$  and a shock-front $\CS$,
which are close to $U^+_0$ and $\mathcal{S}_0$, respectively, such that
\begin{enumerate}
\item[\rm (i)] $U$ satisfies  equations  \eqref{Euler1} in domain $\Om^\s${\rm ;}
\item[\rm (ii)]  The slip condition \eqref{slipcon} holds  along the boundary $\del \mathcal{W}${\rm ;}
\item[\rm (iii)] The Rankine-Hugoniot conditions \eqref{con-RH} as free boundary conditions
hold along the shock-front  $\CS$.
\end{enumerate}
When $U^+_0$ corresponding to a state on arc $\wideparen{TS}$ gives a weak transonic shock,
the problem is denoted by {\rm \textbf{Problem  WT}}, while
the strong transonic shock problem corresponds to arc $\wideparen{TH}$,
denoted by  {\rm \textbf{Problem  ST}}.
\end{problem}

To state our results, we need to introduce the weighed H\"older norms for our subsonic domain $E$,
where $E$ is either a truncated triangular domain or an unbounded domain with the vertex at
origin $\O$ and one side as the wedge boundary.
There are two weights: One is the distance function to origin $\O$, and the other is to
the wedge boundary $\partial \mathcal{W}$.
For any $\xx, \xx'\in  E$,  define
\begin{align*}
&\gd^{\rm o}_\xx := \min(|\xx|,1), &\gd^{\rm o}_{\xx,\xx'} := \min (\gd^{\rm o}_\xx, \gd^{\rm o}_{\xx'}), \\[1mm]
&\gd^{\rm w}_\xx := \min(\textrm{dist}(\xx,\partial \mathcal{W}),1),
  & \gd^{\rm w}_{\xx,\xx'}:= \min (\gd_\xx^{\rm w},\gd_{\xx'}^{\rm w}), \\[1mm]
&\gD_\xx := |\xx|+1, &\gD_{\xx, \xx'}:=\min(\gD_\xx, \gD_{\xx'}),\\[1mm]
&\widetilde{\gD}_\xx := \textrm{dist}(\xx,\partial \mathcal{W})+1,
 &\widetilde{\gD}_{\xx, \xx'}:=\min(\widetilde{\gD}_\xx, \widetilde{\gD}_{\xx'}).
\end{align*}
Let $\ga \in (0,1)$,  $\tau,l, \gamma_1, \g_2 \in \R$ with $\g_1 \ge \g_2$,
and $k$ be a nonnegative integer.
Let $\kk = (k_1, k_2)$ be an integer-valued vector, where $k_1, k_2 \ge 0$,
$|\kk|=k_1 +k_2$, and $D^{\kk}= \po_{x_1}^{k_1}\po_{x_2}^{k_2}$.
We define
\begin{align}
&[ f ]_{k,0;(\tau,l);E}^{(\g_1;\O)(\g_2;\partial \mathcal{W})}\nonumber\\
&\quad :=\sup_{\begin{subarray}{c}
\xx\in E\\
 |\kk|=k
\end{subarray}}
\begin{array}{l}
\big\{(\gd^{\rm o}_{\xx})^{\max\{\g_1 + \min\{k,-\g_2\},0\}}(\gd^{\rm w}_{\xx})^{\max\{k+\g_2,0\}}
\,\gD_\xx^{\tau} \widetilde{\gD}_\xx^{l+k} |D^\kk f(\xx)|\big\},
\end{array}\label{def-normk0}\\[2mm]
&{[ f ]}_{k,\ga;(\tau,l);E}^{(\g_1;\O)(\g_2;\partial \mathcal{W})}
:=\sup_{
\begin{subarray}{c}
   \xx, \xx'\in E\\
 \xx \ne \xx'\\
|\kk|=k
 \end{subarray}}
\Bigg\{\begin{array}{l}
  (\gd^{\rm o}_{\xx,\xx'})^{\max\{\g_1 + \min\{k+\ga,-\g_2\},0\}}(\gd^{\rm w}_{\xx,\xx'})^{\max\{k+\ga+\g_2,0\}}\\[1mm]
\,\,\,\times \gD_{\xx,\xx'}^{\tau} \widetilde{\gD}_{\xx,\xx'}^{l+k+\ga}\frac{|D^\kk f(\xx)-D^\kk
f(\xx')|}{|\xx-\xx'|^\ga}
\end{array}\Bigg\},
\label{def-normkga}
\\
& \|f\|_{k,\ga;(\tau,l); E}^{(\g_1;\O)(\g_2;\partial \mathcal{W})}:= \sum_{i=0}^k {[ f
]}_{i,0 ;(\tau,l);E}^{(\g_1;\O)(\g_2;\partial \mathcal{W})} + {[ f
]}_{k,\ga;(\tau,l);E}^{(\g_1;\O)(\g_2;\partial \mathcal{W})}. \label{def-norm}
\end{align}

For a vector-valued function $\ff=(f_1, f_2, \cdots,f_n )$, we
define
\[
 \|\ff\|_{k,\ga;(\tau,l); E}^{(\g_1;\O)(\g_2;\partial \mathcal{W})} :=
\sum_{i=1}^n  \|f_i\|_{k,\ga;(\tau,l); E}^{(\g_1;\O)(\g_2;\partial \mathcal{W})}.
\]
 Let
\begin{equation}\label{def-C}
C^{k,\ga;(\tau,l)}_{(\g_1;\O)(\g_2;\partial \mathcal{W})}(E)
:= \{ f: \|f\|_{k,\ga;(\tau,l); E}^{(\g_1;\O)(\g_2;\partial \mathcal{W})} <\infty \}.
\end{equation}

\begin{remark}
The requirement that $\g_1 \ge \g_2$ in the definition above means that the regularity up
to the wedge boundary is no worse than the regularity up to the vertex.
When $\g_1 = \g_2$, the $\d^{\rm o}$--terms disappear
so that $(\g_1;\O)$ in the superscript or subscript can be dropped.

If there is no weight $(\g_2;\pw)$ in the superscript, the $\d$--terms for the weights
should be understood as $(\d_{\xx}^{\rm o})^{\max\{k+\g_1,0\}}$
and $(\d_{\xx}^{\rm o})^{\max\{k + \ga +\g_1,0\}}$
in \eqref{def-normk0} and \eqref{def-normkga}, respectively.
When no weight appears in the superscripts of the seminorms
in \eqref{def-normk0}--\eqref{def-normkga},
it means that neither $\delta^{\rm o}$ nor $\delta^{\rm w}$ is present.

For a function of one variable defined on $(0,\infty)$,
the weighted norm $\|f\|^{(\g_2;0)}_{k,\ga;(l);\R^+}$
is understood in the same sense as the definition
above with weight to $\{0\}$ and the decay
at infinity.
\end{remark}

Since the components of $U$
are expected to have different regularity,
we distinguish these variables by defining
${U_1}=(u_1,\rho)$ and $U_2=(w, p)$, where $w=\frac{u_2}{u_1}$.
Let $U_{10}^+$ and $U_{20}^+$ be
the corresponding background subsonic states.

\begin{theorem}[Main Theorem] \label{thm-main}
There are positive constants $\ga , \b,  C_0$, and $\varepsilon$, depending only on the background
states $(U^-_0, U^+_0)$,  such that

\begin{enumerate}
\item[\rm (i)] When $U^+_0\in \wideparen{TS}$, then, for every upstream flow $U^-$ and wedge boundary $x_2=b(x_1)$
satisfying
\[
\| U^- -U^-_0 \|_{2,\ga;(1+\b,0);\Om^-} +\| b' \|^{(-\ga;0)}_{1,\ga; (1+\b);\R^+} \le \ve,
\]
there exists a solution $(U,\s)$ of {\rm \textbf{Problem WT}} satisfying
\begin{equation}
\begin{split}
& \|U - U_0^+ \|_{X } + \| \s' - k_0\|^{(-\ga;0)}_{2,\ga;(1+\gb);\R^+}\\
&\le{}  C_0 \left(\| U^- -U^-_0 \|_{2,\ga;(1+\b,0);\Om^-}
     +\|b'\|^{(-\ga;0)}_{1,\ga; (1+\b);\R^+} \right),
\label{est-U-small-pert}
\end{split}
\end{equation}
where
\begin{equation*}
\|U - U_0^+ \|_{X }:= \|U_1 - U_{10}^+\|^{(-\ga;\del  \mathcal{W})}_{2,\ga;(0,1+\gb);\Om^\s}
+\|U_2 - U_{20}^+\|^{(-\ga;  \O)(-1-\ga;\del  \mathcal{W})}_{2,\ga;(1+\gb,0);\Om^\s};
\end{equation*}

\item[\rm (ii)] When $U^+_0\in \wideparen{TH}$, then, for every upstream flow $U^-$ and wedge boundary $x_2=b(x_1)$ satisfying
\[
\|U^- -U^-_0\|_{2,\ga;(\b,0);\Om^-} + \|b'\|^{(-1-\ga;0)}_{2,\ga; (\b);\R^+}\le \ve,
\]
there exists a solution $(U,\s)$ of {\rm \textbf{Problem ST}} satisfying
\begin{equation}
\begin{split}
& \|U - U_0^+ \|_{X'}
+\|\s' - k_0\|^{(-1-\ga;0)}_{2,\ga;(\gb);\R^+}
\\
&\le{} C_0  \left(\| U^- -U^-_0\|_{2,\ga;(\b);\Om^-}
  +\|b' \|^{(-1-\ga;0)}_{2,\ga; (\b);\R^+}  \right), \label{est-U-small-pert2}
\end{split}
\end{equation}
where
\begin{eqnarray*}
	\|U - U_0^+ \|_{X' }
	:= \|U_1 - U_{10}^+\|^{(-1-\ga;\pw)}_{2,\ga;(0,\gb);\Om^\s}
	+\|U_2 - U_{20}^+\|^{(-1-\ga; \O)}_{2,\ga;(\gb,0);\Om^\s}.
\end{eqnarray*}
\end{enumerate}
The solution $(U,\s)$ is unique within the class such that
  the left-hand side of
\eqref{est-U-small-pert} for {\rm \textbf{Problem WT}} or \eqref{est-U-small-pert2}
for {\rm \textbf{Problem ST}} is less than $C_0 \vae$.
\end{theorem}

\begin{remark}\label{rem:2.2a}
The dependence of constants $\ga , \b,  C_0$, and $\varepsilon$ in Theorem \ref{thm-main}
is described as follows:  $\ga$ and $\b$ depend on $U^-_0$ and $U^+_0$,
but are independent of $C_0$ and $\varepsilon$;
$C_0$ depends on $U^-_0, U^+_0,\ga$, and $\beta$, but is independent of $\ve$;
and $\ve$ depends on all $U^-_0, U^+_0, \ga , \b$, and $C_0$.
\end{remark}

\begin{remark}
The difference in the results of the two problems is that the solution of
{\rm \textbf{Problem WT}} has less regularity at corner ${\bf O}$
and decays faster with respect to $|\xx|$ (or the distance from the wedge boundary)
than the solution of {\rm \textbf{Problem ST}}.
\end{remark}

\begin{remark}
The asymptotic behavior of the subsonic solution can be stated more clearly
in Lagrangian coordinates.
Thus we leave it in the statement of Theorem {\rm \ref{thm-lag}} and Remark {\rm \ref{rk-asymp}}.
\end{remark}

\section{The Problem in Lagrangian Coordinates} \label{sec-lagrange}

From the first equation in \eqref{Euler1}, there exists a unique
stream function $\psi$ in region $\Om^- \cup \Om^{\s}$ such that
$$
\psi_{x_1}= -\rho u_2, \qquad \psi_{x_2} = \rho u_1
$$
with $\psi(\mathbf{0})=0$.
To
simplify the analysis, we employ the following Lagrangian coordinate
transformation:
\begin{equation}\label{def-coord}
\begin{cases}
y_1=x_1,\\[1mm]
y_2= \psi(x_1, x_2),
\end{cases}
\end{equation}
under which the original curved streamlines become straight.
In the new coordinates $\yy=(y_1, y_2)$, we still denote the unknown
variables $U(\xx(\yy))$ by $U(\yy)$ for notational simplicity.

The Euler equations in \eqref{Euler1} in Lagrangian coordinates
become the following equations in divergence form:
\begin{align}\label{eqn-euler1}
&\Big(\frac{1}{\rho u_1} \Big)_{y_1}
    -\Big(\frac{u_2}{u_1} \Big)_{y_2}= 0,\\[1mm]
&\Big(u_1 + \frac{p}{\rho u_1}\Big)_{y_1}
   -\Big(\frac{p u_2}{u_1}\Big)_{y_2}= 0, \label{eqn-euler2}\\[1mm]
&(u_2)_{y_1} + p_{y_2}= 0,\label{eqn-euler3} \\[1mm]
& \Big(\frac{1}{2}|\uu|^2 +\frac{\gc p}{(\gc-1)\rho}\Big)_{y_1}=0.
  \label{eqn-euler4}
\end{align}

Let $\mathcal{T}:=\{y_1= \hs(y_2)\}$ be a shock-front in the $\mathbf{y}$--coordinates.
Then, from the
equations above, we can derive the Rankine-Hugoniot conditions
along $\mathcal{T}$:
\begin{align}\label{con-RH1}
&\Big[\frac{1}{\rho u_1}\Big]=-\Big[\frac{u_2}{u_1} \Big] \hs'(y_2),\\[1mm]
& \Big[ u_1 + \frac{p}{\rho u_1}\Big]=-\Big[\frac{p u_2}{u_1}\Big]
 \hs'(y_2),
  \label{con-RH2}\\[1mm]
&[\,u_2 \,]= [\,p\,] \hs'(y_2),
\label{con-RH3}\\[1mm]
&\Big[\frac{1}{2}|\uu|^2 + \frac{\gc p}{(\gc-1)\rho}\Big]= 0.
\label{con-RH4}
\end{align}

The background shock-front now is $\mathcal{T}_0:=\{y_1=  \hs_0(y_2):= k_1 y_2\}$,
where $k_1 =\frac{k_0}{\rho^+_0u^+_{10}}$.
Without loss of generality, we assume that the supersonic
solution $U^-$ exists in region $\D^-$ defined by
\begin{equation}
\label{def-super-D}
  \D^-  :=   \Big\{\yy: 0<y_1 < \frac{4}{3}k_1y_2 \Big\}.
\end{equation}
Let
\begin{eqnarray}
  \D &=&  \left\{\yy \,:\,  0<  k_1y_2  <y_1   \right\},  \label{def-D}\\
    \mathcal{L}_1&=&  \left\{\yy\,:\, y_1>0, y_2 =0 \right\}, \\
     \mathcal{L}_2&=&  \left\{\yy\,:\, y_1>0,  y_1=  k_1y_2   \right\}\label{def-L2} .
\end{eqnarray}

For a given shock function $\hs(y_2)$, let
\begin{equation}\label{}
\D^{\hs}  = \left\{\yy\,: \, y_2 >0, \hs(y_2)< y_1 \right\}.
\end{equation}

\begin{figure}
 \centering
\includegraphics[height=58mm]{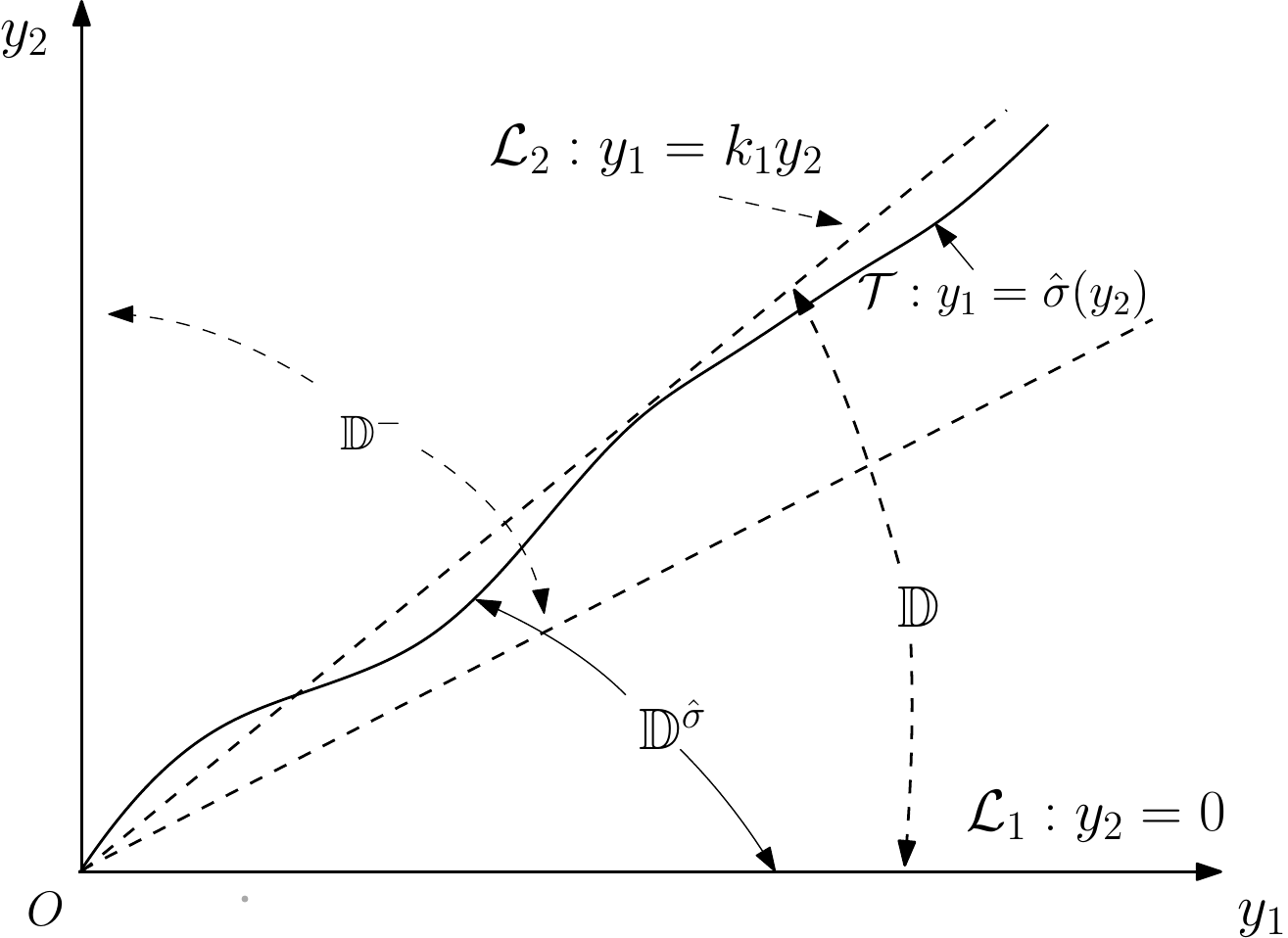}
\caption{Domains $\D$ and $\D^{\hs}$ in Lagrangian coordinates}
\label{new-2}
\end{figure}

Then Theorem \ref{thm-main} can be stated in Lagrangian coordinates as follows:

\begin{theorem} \label{thm-lag}
There exist positive constants $\ga , \b,  C_0$, and $\varepsilon$, depending only on
the background states $U^-_0$ and $U^+_0$, such that,
if the upstream flow $U^-$ for \eqref{eqn-euler1}--\eqref{eqn-euler4} and the wedge boundary function $b(y_1)$
satisfy

\smallskip
\begin{enumerate}
\item[\rm (i)]  $\|U^- -U^-_0\|_{2,\ga;(1+\b,0);\D^-} +\|b'\|^{(-\ga;0)}_{1,\ga; (1+\b);\R^+}\le\ve$ $\,\,$ for {\rm \textbf{Problem WT}};

\smallskip
\item[\rm  (ii)] $\|U^- -U^-_0\|_{2,\ga;(\b,0);\D^-} +\|b'\|^{(-1-\ga;0)}_{2,\ga; (\b);\R^+} \le \ve$ $\,\,$ for {\rm \textbf{Problem ST}},
\end{enumerate}

\smallskip\noindent
then there exist a transonic shock
$\mathcal{T}:=\{y_1= \hs(y_2)\}$ and a subsonic solution $U$ of the Euler
equations \eqref{eqn-euler1}--\eqref{eqn-euler4} satisfying the Rankine-Hugoniot conditions
\eqref{con-RH1}--\eqref{con-RH4} along $\mathcal{T}$ and
the  slip condition $w|_{\mathcal{L}_1}= b'$,
and there exists a limit function $U^\infty(y_2)=(u_1^\infty(y_2), 0, p^+_0, \rho^\infty(y_2))$
and $U_1^\infty(y_2)=(u_1^\infty(y_2), \rho^\infty(y_2))$ such
that  $U$ satisfies the following estimates{\rm :}

\begin{enumerate}
\item[(i)] For {\rm \textbf{Problem WT}},
\begin{equation}
\begin{split}
&\|U - U^\infty \|_{Y }
 +\|\hs' -k_1\|^{(-\ga;0)}_{2,\ga;(1+\gb);\R^+}
 + \|U^{\infty}_1 - U^+_{10}\|_{2,\ga;(1+\gb); \R^+}^{(-\ga;0)}
\\
&\,\,\le{} C_0 \left(\|U^- -U^-_0\|_{2,\ga;(1+\b,0);\D^-}+\|b'\|^{(-\ga;0)}_{1,\ga; (1+\b);\R^+}\right),
 \label{est-U-small-pert3}
 \end{split}
\end{equation}
where
\begin{equation*}
  \|U - U^\infty \|_{Y} := \|U_1 - U_1^\infty \|^{(-\ga;\mathcal{L}_1)}_{2,\ga;(1+\gb,0);\D^{\hs}}
  +\|U_2 - U_{20}^+\|^{(-\ga;  \O)(-1-\ga;\mathcal{L}_1)}_{2,\ga;(1+\gb,0);\D^{\hs}};
 \end{equation*}

\item[(ii)] For {\rm \textbf{Problem ST}},
\begin{equation}
\begin{split}
 &\|U- U^{\infty} \|_{Y' }
 +\|\hs' -k_1 \|^{(-1-\ga;0)}_{2,\ga;(\gb);\R^+} +  \|U^{\infty}_1 - U^+_{10}\|_{2,\ga;(\gb); \R^+}^{(-1-\ga;0)}
\\
 &\,\,\le{} C_0  \left(\|U^- -U^-_0 \|_{2,\ga;(\b,0);\D^-}+\|b'\|_{1,\ga; (\b);\R^+} \right),
 \label{est-U-small-pert4}
 \end{split}
\end{equation}
where
\begin{equation*}
  \|U - U^{\infty} \|_{Y' }:= \|U_1 - U_1^{\infty}\|^{(-1-\ga;\mathcal{L}_1)}_{2,\ga;(\gb,0);\D^{\hs}}
 +\|U_2 - U_{20}^+\|^{(-1-\ga;\O)}_{2,\ga;(\gb,0);\D^{\hs}}.
 \end{equation*}
\end{enumerate}
Moreover, solution $U$ is unique in the class such that the left-hand side of
estimate \eqref{est-U-small-pert3} {\rm (}for {\rm \textbf{Problem WT}}{\rm )}
or \eqref{est-U-small-pert4} {\rm (}for {\rm \textbf{Problem ST}}{\rm )}
is less than $C_0 \ve$. See also Fig. {\rm \ref{new-2}}.
\end{theorem}

\begin{remark}\label{rk-asymp}
In general, the asymptotic limit $U_1^\infty$ is not a constant,
which indicates that $(u_1, \rho)$ does not converge to the background state $(u_{10}^+, \rho_0^+)$
as $y_1 \to  \infty$ (along the streamlines);
while $(u_1, \rho)$ converges to the background state as $y_2 \to \infty$ (transversal to the streamlines away from the wedge).
Such an asymptotic behavior is owing to the conservation of the Bernoulli quantity and the entropy function
along the streamlines, which is different
from that for potential flows.
 \end{remark}

\begin{remark}
Estimates \eqref{est-U-small-pert3}--\eqref{est-U-small-pert4} in Theorem \ref{thm-lag},
together with the Rankine-Hugoniot conditions \eqref{con-RH1},
imply that the coordinate transformation \eqref{def-coord}
is bi-Lipschitz across the shock-front and has the corresponding regularity
in each supersonic or subsonic domain.
Therefore,  Theorem \ref{thm-lag} implies  Theorem \ref{thm-main}.
\end{remark}

\section{Decomposition of the Euler system} \label{sec-decompose}
We now use the left eigenvectors to decompose the Euler
equations \eqref{eqn-euler1}--\eqref{eqn-euler4} into
an elliptic system and two algebraic equations.

Rewrite system  \eqref{eqn-euler1}--\eqref{eqn-euler4} into the
following nondivergence form for $U=(\uu, p,\rho)^\top$:
\begin{equation}
 A(U) U_{y_1} + B(U) U_{y_2} =0, \label{euler-nondiv}
 \end{equation}
 where
\begin{align*}
A(U) &=
\begin{bmatrix}
-\frac{1}{\rho u_1^2} & 0 &0& -\frac{1}{\rho^2 u_1}\\[1.5mm]
1 -  \frac{ p}{ \rho u_1^2} & 0 &  \frac{1}{\rho u_1} & -  \frac{ p}{ \rho^2 u_1}\\[1.5mm]
0  & 1& 0 & 0\\[1mm]
u_1 & u_2 & \frac{\g}{(\g - 1) \rho} & - \frac{\g p}{(\g - 1) \rho^2}
\end{bmatrix},\\[1mm]
B(U) &=
 \begin{bmatrix}
\frac{u_2}{ u_1^2} & -\frac{1}{u_1} &0&0\\[1.5mm]
  \frac{ pu_2}{ u_1^2} & - \frac{p}{u_1} & - \frac{u_2}{ u_1} & 0 \\[1.5mm]
0  & 0& 1& 0\\[1mm]
0& 0 &0 & 0
\end{bmatrix}.
 \end{align*}

Solving $\det (\gl A -B)=0$ for $\gl$, we obtain four eigenvalues:
 \begin{eqnarray*}
 && \l_1 =\l_2 =0,\\
 && \l_{3,4} \equiv \l_{\pm} =-\frac{c\rho}{c^2-u_1^2}\big(c u_2 \mp u_1 \sqrt{c^2 -q^2} i\big),
 \end{eqnarray*}
 where $q= \sqrt{u_1^2 + u_2^2} < c$ in the subsonic region.
 The corresponding left-eigenvectors are
 \begin{eqnarray*}
l_1 & =& (0,0,0,1),\\
l_2 & =&(-p u_1, u_1, u_2, -1),\\
l_{3,4} &=&((\frac{\g p^2}{(\g -1) \rho u_1} - \frac{p u_1}{ \g -1 })\l_{3,4} + \frac{\g p^2 u_2}{(\g - 1) u_1},\\
 &&\,\,-(u_1+\frac{\g p}{(\g -1) \rho u_1})\l_{3,4}-\frac{\g p u_2}{(\g -1)u_1}, \frac{\g p}{\g -1} - u_2 \l_{3,4},  \l_{3,4}).
\end{eqnarray*}

Then
\begin{enumerate}
\item[(i)]
Multiplying equations \eqref{euler-nondiv} from the left by $l_1$ leads to
the same equation \eqref{eqn-euler4}.
This, together with the Rankine-Hugoniot condition \eqref{con-RH4}, implies the Bernoulli law:
\begin{equation}\label{eqn-bernoulli}
\frac{1}{2}q^2 + \frac{\gc p}{(\gc-1)\rho} = B(y_2)
\end{equation}
in both supersonic and subsonic domains, and across the shock-front.
Therefore, $B(y_2)$ can be computed from the upstream flow
$U^-$.  If $u_1$ is a small perturbation of $u_{10}^+$, then $u_1>0$. Therefore,  we can solve \eqref{eqn-bernoulli} for $u_1$:
\begin{equation}\label{eqn-u1}
u_1 = \frac{ \sqrt{2B - \frac{2\g p }{ (\g -1 )\rho }}}{\sqrt{1+ w^2}}
\end{equation}
with $w:=\frac{u_2}{u_1}$.

\item[(ii)] Multiplying system \eqref{euler-nondiv} from the left by $l_2$ gives
\begin{equation} \label{eqn-prho}
\Big(\frac{p}{\rho^\g}\Big)_{y_1} =0.
\end{equation}

\item[(iii)] Multiplying  equations \eqref{euler-nondiv} from the left by $l_3$
and separating the real and imaginary parts of the equation lead to the elliptic system:
\begin{eqnarray}
&&D_R w + e D_I p  =  0,   \label{eqn-wp1}\\
&&D_I w - e D_R p  = 0, \label{eqn-wp2}
\end{eqnarray}
where $e= \frac{\sqrt{c^2-q^2}}{c\rho u_1^2}$ and
$$
\qquad  D_R = \del_{y_1} + \l_R \del_{y_2}, \,\,\,  D_I = \l_I\del_{y_2},\quad
\l_R =-\frac{c^2 \rho u_2}{c^2 - u_1^2 }, \,\,\,\l_I =  \frac{c\rho u_1 \sqrt{c^2 -q^2}}{c^2 - u_1^2}.
$$
\end{enumerate}
Therefore, equations \eqref{eqn-euler1}--\eqref{eqn-euler4} are decomposed into \eqref{eqn-u1}--\eqref{eqn-wp2}.

\smallskip
We will follow the steps below to solve this problem:

\smallskip
1. Given a shock-front $\hs$, introduce a linear system \eqref{eqn-wp-z1}--\eqref{eqn-wp-z2} for iterations;

\smallskip
2. For a given $U$, find $\tU$ by solving the linear system \eqref{eqn-wp-z1}--\eqref{eqn-wp-z2}
with equations \eqref{eqn-u1}--\eqref{eqn-prho} and the corresponding
boundary conditions;

\smallskip
3. Use solution $\tU$ to update the shock-front and obtain $\tilde{\s}$,
so that we construct a map $\mathcal{Q}$ from $(\d U, {\d \hs}^{\prime})$ to $(\d \tU, {\d \tilde{\s}}^{\prime})$;

\smallskip
4. Prove the existence of the solution as a fixed point of $\mathcal{Q}$ by applying
the Schauder fixed point theorem.

\section{Linear Boundary Value Problem for Iterations} \label{sec-linearize}
For a given shock-front $\hs$, the subsonic domain $\D^{\hs}$ depends on $\hs$.
For the convenience of solving the problem, we make the following coordinate transformation
to change the domain from $\D^{\hs}$ to $\D$:
\begin{equation}\label{eqn-tranz}
\left \{ \begin{aligned}
z_1 &= y_1 - \d \hs(y_2),\\
z_2&= y_2,
\end{aligned} \right.
\end{equation}
where $\d\hs(y_2)=\hs(y_2)-\hs_0(y_2)$.
In the $\zz$--coordinates, $U(\yy)$ becomes $U_{\hs}(\zz)$, depending on $\hs$.
When there is no ambiguity, we may omit the subscript and still denote $U_{\hs}(\zz)$ by $U(\zz)$.
However, the upstream flow $U^-$ involves an unknown variable explicitly depending on $\hs$:
$$
U^-_{\hs}(\zz) = U^-(z_1+ \d \hs(z_2), z_2),
$$
where $U^-$ is the given upstream flow in the $\yy$--coordinates.
Hence, equations \eqref{eqn-wp1}--\eqref{eqn-wp2} become the following equations
in  the $\zz$-coordinates:
 \begin{eqnarray}
&&\widetilde{D}_R w + e \widetilde{D}_I p = 0,   \label{eqn-wp-z1}\\
&&\widetilde{D}_I w - e \widetilde{D}_R p = 0,\label{eqn-wp-z2}
\end{eqnarray}
 where
 \begin{align*}
 \widetilde{D}_R= (1-\gd \hs' \gl_R)\partial_{z_1} + \gl_R\partial_{z_2}, \qquad
  \widetilde{D}_I= \gl_I (-\gd \hs' \partial_{z_1} +  \partial_{z_2}).
 \end{align*}
 Using system \eqref{eqn-wp-z1}--\eqref{eqn-wp-z2} to solve for $(p_{z_1}, p_{z_2})$ yields
 the linear system for iterations:
 \begin{align}
(\d \tilde{p})_{z_1} & =\frac{\gl_R-\d \hs'(\gl_R^2 +\l_I^2)}{e \gl_I }(\d \tilde{w})_{z_1}+\frac{\gl_R^2 +\l_I^2}{e \gl_I }(\d \tilde{w} )_{z_2},
   \label{eqn-wp-lin1}\\
 (\d \tilde{p})_{z_2} & = -\frac{(1-\d \hs'\gl_R)^2 +(\d \hs' \l_I)^2}{e \gl_I }(\d \tilde{w})_{z_1}
  -\frac{\gl_R - \d \hs'(\gl_R^2 + \l_I^2)}{e \gl_I }(\d \tilde{w})_{z_2}.\label{eqn-wp-lin2}
\end{align}

In the $\zz$-coordinates, the Rankine-Hugoniot conditions \eqref{con-RH1}--\eqref{con-RH4} keep the same form,
except that $\hs'(y_2)$ is replaced by $\hs'(z_2)$ and $U^-$ is replaced by $U^-_{\hs}$ along line $\mathcal{L}_2$.
Among the four Rankine-Hugoniot conditions, \eqref{con-RH4} is used in the Bernoulli law.
From condition \eqref{con-RH3}, we have
\begin{equation} \label{eqn-sigmaprime}
 \hs' (z_2) = \frac{[ u_1 w ]}{[p]}(k_1 z_2, z_2),
\end{equation}
which will be used to update the shock-front later.
Now, because of \eqref{eqn-u1}, we can use $\bar{U}=(w, p, \rho)$ as the unknown variables along $\mathcal{L}_2$.
Using \eqref{eqn-sigmaprime} to eliminate $\hs'$ in conditions  \eqref{con-RH1}--\eqref{con-RH2} gives
\begin{eqnarray}
&&G_1(U_{\hs}^-, \bar{U}) := [p] \Big[ \frac{1}{\rho u_1 }\Big] + [w] [u_1 w] = 0, \label{con-G1}\\
&&G_2(U_{\hs}^-, \bar{U}) := [p] \Big[ u_1 + \frac{p}{\rho u_1 }\Big] + [pw] [u_1 w] = 0.\label{con-G2}
\end{eqnarray}
We linearize the conditions above as
\begin{equation}
\grad_{\bar{U}} G_i(U_0^-, \bar{U}_0^+)\cdot \d \widetilde{\bar{U}}
= \grad_{\bar{U}} G_i(U_0^-, \bar{U}_0^+)\cdot \d \bar{U} - G_i(U_{\hs}^-, \bar{U}),
\label{con-lin-gi}
\end{equation}
denoted by
\begin{equation}
b_{i1}\d \tilde{w}  + b_{i2}\d \tilde{p}  +b_{i3}\d \tilde{\rho}  =g_i (U_{\hs}^-, \bar{U}), \qquad i=1,2,\label{con-gi}
\end{equation}
where
\begin{eqnarray}
&&(b_{i1},b_{i2},b_{i3}) := \grad_{\bar{U}} G_i(U_0^-, \bar{U}_0^+),\\
&&g_i (U_{\hs}^-, \bar{U}):= \grad_{\bar{U}} G_i(U_0^-, \bar{U}_0^+)\cdot \d \bar{U}
- G_i(U_{\hs}^-, \bar{U}). \label{eqn-gi}
\end{eqnarray}
 Using the two conditions \eqref{con-gi}, for $i=1,2$, to eliminate $\d \tilde{\rho} $ leads to
 \begin{equation} \label{con-g30}
 (b_{11}b_{23} - b_{21} b_{13}) \d \tilde{w}  + (b_{12}b_{23} - b_{22}b_{13})\d \tilde{p}   =b_{23}g_1 -b_{13}g_2.
 \end{equation}
 A direct calculation shows
 \begin{eqnarray*}
&&b_{11}b_{23} - b_{21} b_{13} \\
&&=  (-u_{20}^-)[p_0] \Big( \frac{\g p_0^+}{(\g -1)(\rho_0^+)^2 u_{10}^+} + \frac{p_0^-}{u_{10}^-}\big(\frac{1}{(\rho_0^+)^2 }
 +  \frac{\g p_0^+}{(\g -1)(\rho_0^+)^3 (u_{10}^+)^2}\big)\Big)\\
&& > 0.
\end{eqnarray*}
 Therefore, condition \eqref{con-g30} becomes
  \begin{equation} \label{con-g3}
 \d \tilde{w} + b_1\d \tilde{p}=g_3,
 \end{equation}
 where
 \begin{equation}\label{5.14a}
 b_1=\frac{b_{12}b_{23}-b_{22}b_{13}}{b_{11}b_{23} - b_{21} b_{13}},  \qquad
 g_3 = \frac{b_{23}g_1-b_{13} g_2}{b_{11}b_{23} - b_{21} b_{13}}.
 \end{equation}
 \begin{figure}
 \centering
\includegraphics[height=65mm]{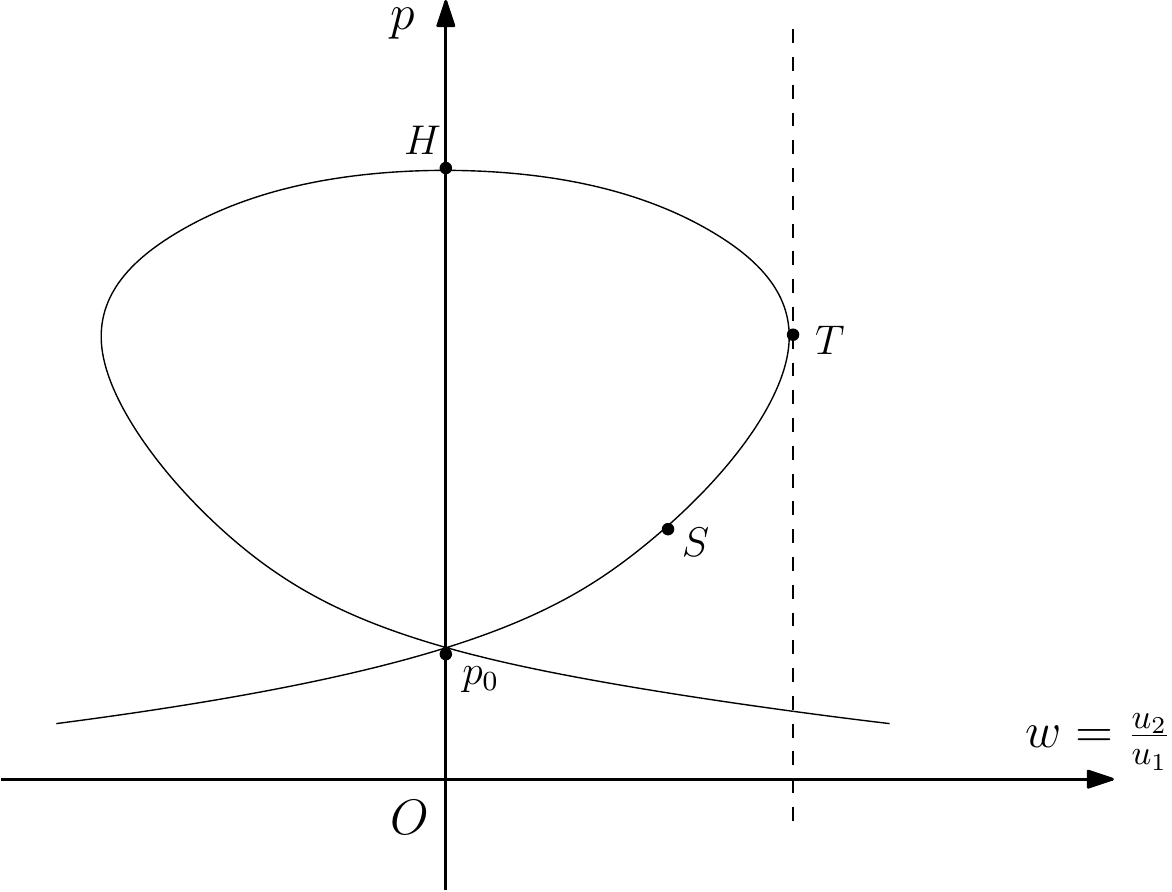}
\caption{The shock polar in the $(w, p)$--variables}
\label{wppolar}
\end{figure}

\begin{remark}\label{signb1}
The shock polar is a one-parameter curve determined by the Rankine-Hugoniot conditions.
If $p$ is used as the parameter, by equation \eqref{con-g3},
we obtain that $\d w= -b_1 \d p + g_3(\d p)$, which shows that $-b_1 \d p $ is the linear term and $g_3(\d p)$ is the higher order term.
From Fig.~\ref{wppolar}, we know that $w(p)$ is decreasing in $p$  on arc $\wideparen{TH}$
and increasing on $\wideparen{TS}$.
Therefore, it is easy to see that $b_1>0$ corresponds to the state on arc $\wideparen{TH}$,
$b_1<0$ to $\wideparen{TS}$, and $b_1 =0$ at the tangent point $T$.
\end{remark}

We compute
$$
b_{13}= -[p_0]\Big(\frac{p_0^+}{(\rho_0^+)^2 u_{10}^+} +\frac{\g p_0^+}{(\g -1)(\rho_0^+)^3 (u_{10}^+)^3}\Big)< 0.
$$
Thus condition \eqref{con-gi} for $i=1$ can be rewritten as
\begin{equation} \label{con-g4}
 \d \tilde{\rho}  =g_4 - b_2 \d  \tilde{w}   -b_3  \d \tilde{p} ,
\end{equation}
where $g_4 =\frac{g_1}{b_{13}}, b_2= \frac{b_{11}}{b_{13}}$, and $b_3=\frac{b_{12}}{b_{13}}$.

We notice that conditions \eqref{con-g3}--\eqref{con-g4} are equivalent to conditions \eqref{con-gi} for $i=1,2$.

\section{Key elliptic estimates} \label{sec-keylemma}

Consider the elliptic equation
\begin{equation}
(a_{ij} v_{z_i})_{z_j}=0    \qquad \text{ in }   \D,
\label{eqn-v}
\end{equation}
with boundary conditions:
\begin{eqnarray}
&& v|_{\mathcal{L}_1} = g_5(z_1),  \label{con-g5} \\
&& \left.\frac{ Dv}{D \boldsymbol{\nu}}\right|_{\mathcal{L}_2}
 \equiv \nabla v \cdot   \boldsymbol{\nu} |_{\mathcal{L}_2} =  g_6(z_2), \label{con-g6}
\end{eqnarray}
where $\D$ is the unbounded triangular domain with two boundaries $\mathcal{L}_1$ and $\mathcal{L}_2$
defined by \eqref{def-D}--\eqref{def-L2}, and
$\boldsymbol{\nu} =(\nu_1, \nu_2)$ is a constant vector with  $|\boldsymbol{\nu} |=1$.
Let $\omega_0 \in (0, \frac{\pi}{2})$ be the angle between $\mathcal{L}_1$ and $\mathcal{L}_2$,
and let $\nu_n=\boldsymbol{\nu} \cdot (-\sin \omega_0, \cos \omega_0)$
and $\nu_t= \boldsymbol{\nu} \cdot ( \cos \omega_0,\sin \omega_0)$
be the normal and tangent components of $\boldsymbol{\nu}$, respectively.
Note that $(-\sin \omega_0, \cos \omega_0)$ is the outer normal to ${\mathcal L}_2$,
and $(\cos \omega_0, \sin \omega_0)$ is tangent to ${\mathcal L}_2$ directed away
from the corner on domain $\D$.
We assume that $\nu_n >0$.

\begin{lemma} \label{thm-lemma}
Consider the boundary value problem \eqref{eqn-v}--\eqref{con-g6}.

\begin{enumerate}
\item[(i)] When $\nu_t <0$,
there exist suitably small $\ga, \gb \in (0,1)$, depending only on $\boldsymbol{\nu}$ and $\omega_0\in (0,\frac{\pi}{2})$,
such that, if
\begin{align}\label{con-aij}
\| a_{ij} - \d_{ij} \|^{(-\ga;\mathcal{L}_1)}_{1,\ga;(0,1+\gb);\D} \le \d
\end{align}
for a suitably small constant $\d >0$ depending only
on $\boldsymbol{\nu}, \omega_0, \ga$, and $\gb$, and
\[
\d_{ij} = \begin{cases}
1 & \quad \text{ if } \,\,\, i=j,\\
0 & \quad \text{ if } \,\,\,  i\neq j,
\end{cases}
\]
$g_5 \in C_{(-\ga;0)}^{1,\ga; (1+\b)}(\R^+)$, and
$g_6 \in C_{(1-\ga;0)}^{1,\ga; (2+\b)}(\R^+)$,
then there exists a unique solution $v \in  C_{(-\ga;  \O)(-1-\ga;\mathcal{L}_1)}^{2,\ga;(1+\gb,0)}(\D)$
of problem \eqref{eqn-v}--\eqref{con-g6}.
Furthermore, there exists a constant $C >0$, depending only on $\boldsymbol{\nu}, \omega_0, \ga$, and $\gb$,
such that the following estimate holds{\rm :}
\begin{eqnarray}\label{est-v1}
\| v \|^{(-\ga;  \O)(-1-\ga;\mathcal{L}_1)}_{2,\ga;(1+\gb,0);\D}
 &\le& C \big(\|g_5 \|^{(-\ga;0)}_{1,\ga; (1+\b);\R^+}
 + \| g_6\|^{(1-\ga;0)}_{1,\ga; (2+\b);\R^+}\big).
\end{eqnarray}

\item[(ii)] When $\nu_t \ge 0$, there exist suitably small $\ga, \gb \in (0,1)$,
depending only on $\boldsymbol{\nu}$ and $\omega_0\in (0, \frac{\pi}{2})$, such that, if
 \begin{align}\label{con-aij2}
 \| a_{ij} - \d_{ij} \|_{1,\ga;(0,\gb);\D} \le \d
 \end{align}
for  a suitably small constant $\d >0$ depending only on $\boldsymbol{\nu}, \omega_0, \ga$, and $\gb$,
$g_5 \in C_{(-1-\ga;0)}^{2,\ga; (\b)}(\R^+)$, and $ g_6 \in C_{(-\ga;0)}^{1,\ga; (1+\b)}(\R^+)$,
then there exists a unique solution $v \in  C_{(-1-\ga;\O)}^{2,\ga;(\gb)}(\D)$ satisfying the following estimate{\rm :}
\begin{eqnarray}
\| v \|^{(-1-\ga;\O)}_{2,\ga;(\gb);\D}
\le C \big(\|g_5 \|^{(-1-\ga;0)}_{2,\ga; (\b); \R^+}
+ \| g_6 \|^{(-\ga;0)}_{1,\ga; (1+\b);\R^+}\big),
\label{est-v}
\end{eqnarray}
where $C >0$ is a constant, depending only on $\boldsymbol{\nu}, \omega_0, \ga$, and $\gb$.
\end{enumerate}
\end{lemma}

In the following estimates, all constants $C, C_i, c_i$, etc. are generic positive constants
depending only on the background
states $U^-_0$ and $U^+_0$ (or $\boldsymbol{\nu}$
and $\omega_0$ in Lemma \ref{thm-lemma}), $\ga$, and $\gb$.

\subsection{$C^0$--estimates} \label{C0estimates}

We first prove part $\mathrm{ (i)} $ of Lemma \ref{thm-lemma}.

We truncate domain $\D$ by line $L_R=\{\zz\,:\, z_1=R\}$, $R>2k_1$,
into a triangle $\D^R=\{\zz\,:\, 0<  k_1z_2  <z_1<R \}$
and prescribe the following boundary condition:
\begin{equation}\label{con-gR}
v|_{L_R} = g_5(R).
\end{equation}
Since $\D^R$ is a bounded domain, we can start with a Neumann condition on $\mathcal{L}_2$
and Dirichlet conditions on $\mathcal{L}_1$ and $L_R$,
and then use the continuity method to prove that there exists a unique solution
$v_R \in C^0(\overline{\D^R})\bigcap C^{2,\ga}(\D^R)$  ({\it cf}. Theorem 1 in \cite{Lieberman}).
The process is standard, based on the {\it apriori} estimates for $v_R$.
We will focus on obtaining the desired estimates of $v_R$, independent of $R$.

Denote
\[
M:=\|g_5 \|^{(-\ga;0)}_{1,\ga; (1+\b);\R^+} + \| g_6 \|^{(1-\ga;0)}_{1,\ga; (2+\b);\R^+}.
\]
The $C^0$--estimates consist of two parts -- corner estimates and decay estimates.

\medskip
\textbf{Corner estimates}. Let $\bar{v}_R (\zz) := v_R(\zz) - g_5(0)$.
Assume
$M>0$ (otherwise, the maximum principle applied to the zero boundary conditions implies a trivial solution),
and set $\bar{\t}:= (\ga +\tau) \t +\t_0$.
Define a comparison function:
\[
v_1= CM \big(r^{\ga} \sin \bar{\t} +z_2^{\ga}\big),
\]
where $(r,\t)$ are the polar coordinates.
Choose  $\tau, \t_0 >0$ suitably small, so that $(\ga +\tau) \om_0 +\t_0 < \frac{\pi}{2}$.

Now we estimate $(a_{ij}(v_1)_{z_i})_{z_j}$ in the following steps. First,
\begin{align*}	
\Delta v_1 &= C M \big((\ga^2 -(\ga +\tau)^2)r^{-2+\ga}\sin \bar{\t} +  \ga(\ga-1) z_2^{-2+\ga}\big)\\
&\le -Cc_1Mr^{-2+\ga}\big(1+ (\sin \t)^{-2+\ga}\big).
\end{align*}
Condition \eqref{con-aij} implies that
\begin{align*}
|(a_{ij} -\d_{ij})\partial_{ij} v_1|& \le CMC_1\d r^{-2+\ga} (\sin \t)^{-2+\ga}.
\end{align*}
Also,
\begin{align*}
(a_{ij})_{z_j} & = O(\d) \big(\max(z_2, 1)\big)^{-2-\gb} \big(\min(z_2, 1)\big)^{-1+\ga}\\
 & = O(\d) \big(\max(z_2,1)\big)^{-1-\ga- \gb}
  \big( \max(z_2,1) \min(z_2,1)\big)^{-1+\ga}\\
& = O(\d)\big(\max(z_2,1)\big)^{-1-\ga- \gb} z_2^{-1+\ga}\\
& = O(\d)(z_2+1)^{-1-\ga- \gb} r^{\ga-1} (\sin \theta)^{-1+\ga}\\
& = O(\d)r^{-1} (\sin \theta)^{-1} .
\end{align*}
This gives rise to the following estimate:
\begin{align*}
|(a_{ij})_{z_j} (v_1)_{z_i} |
& = O(\d) M r^{-1} (\sin \theta)^{-1}\big(r^{-1+\ga} + z_2^{-1+\ga}\big)\\
& \le C_2 M \d r^{-2+\ga} (\sin \t)^{-2+\ga} .
\end{align*}
The estimate above yields
\begin{align*}
\big(a_{ij}(v_1)_{z_i}\big)_{z_j}
=\big(\Delta
   +(a_{ij} - \d_{ij})\partial^2_{z_i z_j}\big)v_1
   +(a_{ij})_{z_j} (v_1)_{z_i}\le 0,
\end{align*}
if $\d$ is chosen sufficiently small.

On the boundaries, we compute
\begin{align*}	
\left. \frac{Dv_1}{D \boldsymbol{ \nu}} \right|_{\mathcal{L}_2}
 &= C Mr^{-1+\ga} \big(\nu_n(\ga +\tau)\cos \bar{\t} + \nu_t \ga \sin \bar{\t}
  + \ga (\sin \t)^{-1+\ga}(\nu_t\sin \t+ \nu_n \cos \t)\big)\big|_{\t =\omega_0}\\[1mm]
 & \ge  C c_2Mr^{-1+\ga}  \qquad\,\, (\text{by choosing a suitably small } \ga )\\
 & > g_6   \qquad\qquad\qquad  (\text{by choosing a suitably large } C ), \\[2mm]
v_1 |_{\mathcal{L}_1}& = CM z_1^\ga \ge g_5(z_1) -g_5(0),\\[1mm]
v_1 |_{L_R}& \ge  g_5(R) -g_5(0).
\end{align*}

Therefore, by the comparison principle, we conclude
\[
\bar{v}_R \le v_1 .
\]
By adding a negative sign to $v_1$, we obtain that $\bar{v}_R \ge -v_1$. Thus, we have
\begin{align} \label{est-vbarcorner}
|\bar{v}_R(\zz)| &\le CM|\zz|^{\ga} \qquad\qquad\quad \mbox{for any $\zz \in \D^R$},  \\[1mm]
|v_R(\zz)| &\le CM(1+ |\zz|^{\ga} ) \qquad \,\,\,\mbox{for any $\zz \in \D^R$}.\label{est-vcorner}
\end{align}
In particular, for $\zz \in \D^{2 k_1}$, we have
\begin{equation}
|v_R(\zz)| \le CM.
\end{equation}

\medskip
\textbf{Decay estimates}. Now we estimate the decay rate of $v_R$ in $\D^R \backslash \D^{k_1}$.
Denote $\bar{\t}:=(1+\gb +\tau)\t + \t_0$, and
let
\[
v_2 (\zz):= M r^{-1- \gb}\left(C_3 \sin \bar{\t}  +C_4 (\sin \t)^{\ga} \right).
\]
For $\zz \in \D^R \backslash \D^{k_1} $, we calculate
\begin{align*}	
\Delta v_2(\zz) &= M r^{-3- \gb}\big\{C_3 ((1+\gb)^2 -(1+\gb +\tau)^2) \sin\bar{\t} \\
&\qquad\qquad\quad\,\,   +  C_4 \left( -\ga(1-\ga) (\sin \t)^{-2+\ga} + ((1+\gb)^2- \ga^2)(\sin \t)^{\ga} \right) \big\}\\
& \le  -CM r^{-3-\gb}  (\sin \t)^{-2+\ga}
\end{align*}
by adjusting $\frac{C_3}{C_4}$ suitably large. Then
\begin{align*}
|(a_{ij})_{z_j} (v_2)_{z_i} |
 = O(\d)r^{-1} (\sin \theta)^{-1}r^{-2-\gb}  (\sin \theta)^{-1+\ga}
 \le C_5 \d r^{-2-\gb} (\sin \t)^{-2+\ga}
\end{align*}
implies that
\begin{align*}
(a_{ij}(v_2)_{z_i})_{z_j}
& \le 0
\end{align*}
for a sufficiently small $\d$.

Moreover, using $\nu_n>0$, $\nu_t<0$, and $\bar\theta\in (0, \frac{\pi}{2})$,
we have
\begin{align*}	
\left. \frac{Dv_2}{D \boldsymbol{ \nu}} \right|_{\mathcal{L}_2}
&= Mr^{-2- \gb} \big\{ C_3(\nu_n(1+\gb +\tau)\cos\bar{\t} -\nu_t (1+\gb) \sin\bar{\t})  \\
& \qquad\qquad\quad\,\,+ \left.C_4 (\nu_n \ga (\sin \t)^{-1+\ga} \cos \t
     -\nu_t (1+\gb)(\sin \t)^\ga) \big\}\right|_{\t = \om_0}\\[1mm]
&\ge  CMr^{-2-\gb}  \qquad\quad (\text{for large $C_3$ and $C_4$})\\
& >  g_6
\end{align*}
and
$$
v_2 |_{\mathcal{L}_1\cup L_R \cup L_{k_1}} \ge  v_R |_{\mathcal{L}_1\cup L_R \cup L_{k_1}}.
$$
By the comparison principle, we conclude
\begin{equation} \label{est-C0decay}
|v_R(\zz)| \le CM|\zz|^{-1-\gb}   \qquad  \mbox{ for }\ddot{}  \zz \in \D^R \backslash \D^{k_1},
\end{equation}
which yields the following $C^0$--estimate:
\begin{equation}
\|v_R\|_{0,0;(1+\gb,0); \D^R} \le CM.
\end{equation}

\medskip
\subsection{$C^{1,\ga}$--estimates} \label{C1estimates}

Since we will let $R$ approach to $\infty$ eventually, the estimates in $\D^{\frac{R}{2}}$ will be sufficient.

Our estimates are based on the standard Schauder interior or boundary estimates in the discs with appropriate scalings;
{\it cf.} Gilbarg-Trudinger \cite{gt}.
On the other hand, we have different scalings for the corner and away from the corner.

\medskip
{\bf Corner estimates.} First, we focus on
the estimates near corner $\O$.
For any point $\zz^0 \in \D^{k_1} $ with polar coordinates $(r_0, \t_0)$,
we divide the situation into three cases:
$\frac{\omega_0}{4} \le \t_0  \le \frac{3 \omega_0}{4}$,
$\frac{3\omega_0}{4}< \t_0  <\omega_0$,
and
$0 < \t_0 < \frac{\omega_0}{4}$.

\medskip
\textit{Case 1}: $\frac{\omega_0}{4} \le \t_0  \le \frac{3 \omega_0}{4}$.
Let $\bar{r} = \frac{r_0}{4} \sin (\frac{\omega_0}{4})$ and $B_{n\bar{r}}= B_{n\bar{r}}(\zz^0)$
for $n\in \N$.
We rescale $B_{n\bar{r}}$ into $B_{n}:= B_{n}(\O)$ by the coordinate transformation:
\[
\yy = \frac {\zz - \zz^0}{\bar{r}}.
\]
Let $\tilde{v}(\yy) = v_R(\zz^0 + \bar{r}\yy)$ and $\hat{v}(\yy) = \tilde{v}(\yy) - g_5(0)$.
By the $C^{0}$--estimate near the corner, we have
\[
|\hat{v}(\yy)| \le CM \bar{r}^\ga.
\]
Equation \eqref{eqn-v} becomes
\[
(\tilde{a}_{ij}\hat{v}_{y_i})_{y_j} =0,
\]
where $\tilde{a}_{ij}(\yy) = a_{ij}(\zz^0 + \bar{r}\yy )$. Since $\bar{r} \le \sqrt{1+ k_1^2}$, it is easy to see that
\begin{align*}
&[\tilde{a}_{ij}]_{0,\ga; B_2}= \bar{r}^{\ga} [a_{ij}]_{0,\ga; B_{2\bar{r}}} \le C\d,\\[1mm]
&\|\tilde{a}_{ij}\|_{0,\ga; B_2} \le \Lambda, \\[1mm]
&\tilde{a}_{ij}(\yy)\xi_i\xi_j \ge \l |\xi|^2 \qquad \mbox{for $\yy \in B_2$}
\end{align*}
for suitably small $\d$, where $\l, \Lambda >0$ are constants depending only on
$\omega_0, \boldsymbol{\nu},\ga$, and $\gb$.

We apply the Schauder  interior estimate (\textit{cf.} Theorem 8.32 in \cite{gt}) to obtain
\begin{equation*}
\|\hat{v} \|_{1,\ga; B_1} \le C\| \hat{v} \|_{0,0;B_2} \le  CM \bar{r}^\ga.
\end{equation*}

Let $\Omega$ be a domain, let $u$ be a function defined in $\Omega$, and set
$d:= \mbox{diam}\, \Omega$.
We define the following norm $\| \cdot \|'$:
\begin{eqnarray*}
&&\| u\|'_{k; \Omega} = \sum_{j=0}^k d^j[ u]_{j,0;\Omega},\\
&&\| u\|'_{k,\ga;\Omega} =\| u\|'_{k; \Omega}+ d^{k+\ga}[ u]_{k,\ga;\Omega}.
\end{eqnarray*}
Then we obtain the estimate for $\bar{v}_R:= v_R - g_5(0)$:
\begin{equation*}
\|\bar{v}_R  \|'_{1,\ga; B_{\bar{r}}} \le CM \bar{r}^\ga,
\end{equation*}
which implies
\begin{equation}\label{est-in}
\|v_R  \|^{(-\ga;\O)}_{1,\ga; B_{\bar{r}}} \le CM .
\end{equation}

\smallskip
{\it Case 2}: $\frac{3\omega_0}{4}< \t_0  <\omega_0$.
Let $\bar{r}=r_0\sin(\frac{\omega_0}{4}), B_{n \bar{r}}= B_{n\bar{r}}(\zz^0), B_{n\bar{r}}^+= B_{n \bar{r}}\cap \D$,
and $T= B_{2\bar{r}}\cap \mathcal{L}_2$.
We use the same scaling as in {\it Case 1}.
Then the boundary estimates for the Poisson equation with the oblique derivative conditions (see Theorem 6.26 in \cite{gt}),
followed by the technique of freezing the coefficients ({\it cf.} Lemma 6.29 in \cite{gt}), imply that
\begin{equation}\label{est-bdobliq}
\|\bar{v}_R \|'_{1,\ga; B_{\bar{r}}^+}
\le C \big(\| \bar{v}_R \|_{0,0;B_{2\bar{r}}^+} +  \bar{r} \| g_6 \|'_{0,\ga;T} \big)
\le CM \bar{r}^{\ga}.
\end{equation}

\smallskip
{\it Case 3}: $0 < \t_0 < \frac{\omega_0}{4}$. Now  $\bar{r}$ and $B_{n\bar{r}}^+$ are defined in the same fashion
as in {\it Case 2}, while $T= B_{2\bar{r}}\cap \mathcal{L}_1$.
The Schauder boundary estimates for the Dirichlet conditions give rise to the $C^{1,\ga}$--estimates
near boundary $\mathcal{L}_1$ ({\it cf.} Corollary 8.36 in \cite{gt}):
\begin{equation}\label{est-bd1alpha}
\|\bar{v}_R \|'_{1,\ga; B_{\bar{r}}^+} \le C \big(\|\bar{v}_R \|_{0,0;B_{2\bar{r}}^+} + \| g_5 - g_5(0) \|'_{1,\ga;T}\big)
\le CM \bar{r}^{\ga}.
\end{equation}

\smallskip
Therefore, estimates \eqref{est-in}--\eqref{est-bd1alpha} in the cases above give the desired corner estimate:
\begin{equation}\label{est-corner1ga}
\|v_R \|_{1,\ga;\D^{k_1}}^{(-\ga;\O)} \le C M.
\end{equation}

\medskip
{\bf Decay estimates.} Now we consider the domain away from the corner: $\D^*:=\D^{\frac{R}{2}} \backslash \D^{k_1}$.
The estimates below follow the similar way to the corner estimates, but with a different scaling.

For any $\zz^0 = (z^0_1, z^0_2)\in \D^*$,
set $\l_0 =\frac{1}{2\sqrt{1+ k_1^2}}$.
Then we consider two cases: $z^0_2 < \l_0$ and  $z^0_2 \ge \l_0$.

\medskip
 {\it Case 1}: $z^0_2 < \l_0$.
 Set $B_n = B_{n\l_0}(\zz^0)$, $B_n^+ = B_n \cap \D$, and $T= B_{2} \cap \mathcal{L}_1$.
 Similarly,  the Schauder boundary estimate yields
 \begin{equation}\label{est-decay1ga1}
 \|v_R \|'_{1,\ga; B_1^+} \le C \big(\| v_R \|_{0,0;B_2^+} +   \| g_5 \|'_{1,\ga;T}\big)
 \le CM |\zz^0|^{-1-\gb},
 \end{equation}
by using \eqref{est-C0decay}.

 {\it Case 2}: $z^0_2 \ge \l_0$.
 Set $\bar{r} = \frac{z_2^0}{2},  B_{n\bar{r}} = B_{n\bar{r}}(\zz^0)$, $B_{n\bar{r}}^+ = B_{n\bar{r}} \cap \D$,
 and $T= B_{2\bar{r}} \cap \mathcal{L}_2$.
 Similar to the $C^{0}$--estimates away from the corner in \S \ref{C0estimates},
 we rescale to the unit disc by the coordinate transformation $\zz =\zz^0 + \bar{r}\yy$
 and then do either the Schauder boundary or the interior estimates for $\tilde{v}(\yy) = v_R(\zz^0 + \bar{r}\yy)$.
 Since
 \begin{align*}
 [\tilde{a}_{ij}]_{0,\ga;B_2^+}
 = \bar{r}^{\ga} [a_{ij}]_{0,\ga; B^+_{2\bar{r}}}
 \le C\bar{r}^{-1-\gb} [a_{ij}]_{0,\ga;(0,1+\gb); B^+_{2\bar{r}}}
 \le C\d,
 \end{align*}
we obtain the following estimate:
 \[
  \|\tilde{v} \|_{1,\ga; B_1^+} \le C \big(\|\tilde{v}\|_{0,0;B_2^+} + \|\tilde{g}_6\|_{0,\ga;\tilde{T}}\big),
 \]
 where $\tilde{g}_6$ and $\tilde{T} $ are the rescaled function of  $g_6$ and the rescaled boundary of $T$,
 respectively.

Scaling back to $B^+_{\bar{r}}$ leads to
\begin{align}\label{est-decay1ga2}
\|v_R \|'_{1,\ga; B_{\bar{r}}^+}& \le C \big(\|v_R \|_{0,0;B_{2\bar{r}}^+} + \bar{r} \| g_6 \|'_{0,\ga;T}\big)
\le CM |\zz^0|^{-1-\gb}.
\end{align}
Estimates \eqref{est-decay1ga1}--\eqref{est-decay1ga2} give rise to the $C^{1,\ga}$--estimate in $\D^*$:
\begin{equation}\label{est-D*1ga}
\|v_R \|_{1,\ga;(1+\gb,0);\D^*} \le C M .
\end{equation}

Combining estimates \eqref{est-corner1ga} in $\D^{k_1}$ with estimate \eqref{est-D*1ga} in $\D^*$
renders the following $C^{1,\ga}$--estimate in
 $\D^{\frac{R}{2}}$:
 \begin{equation}\label{est-DR1alpha}
\|v_R \|_{1,\ga;(1+\gb,0);\D^{\frac{R}{2}}}^{(-\ga; \O)} \le C M .
\end{equation}

\smallskip
\subsection{$C^{2,\ga}$--estimates}
For the $C^{2,\ga}$--estimates, we rewrite equation \eqref{eqn-v} into the following non-divergence form:
\begin{equation}\label{eqn-v2ga}
a_{ij}v_{z_iz_j} + (a_{ij})_{z_j} v_{z_i} = 0.
\end{equation}
Following the same argument as the $C^{1,\ga}$--estimates in \S \ref{C1estimates},
we let $\bar{r}=\frac{z_2^0}{4},  B_{n\bar{r}}=B_{n\bar{r}}(\zz^0)$, and $B_{n\bar{r}}^+=B_{n\bar{r}}\cap\D$.
For $z_2^0 \ge 1$, we follow the same procedure as in the $C^{1,\ga}$--estimates to conclude
\begin{align}\label{est-decay2ga2}
\|v_R \|'_{2,\ga; B_{\bar{r}}^+}
& \le C \big(\|v_R \|_{0,0;B_{2\bar{r}}^+} + \bar{r}\|g_6 \|'_{1,\ga;T}\big)
\le CM |\zz^0|^{-1-\gb}.
\end{align}
To obtain the estimates for $z_2^0 <1$, we set
\begin{align*}
&\bar{\zz}^0 = (z_1^0 , \frac{z^0_2}{4}),\\
&\bar{v}(\zz) = v_R(\zz) - v_R(\bar{\zz}^0) - \nabla v_R(\bar{\zz}^0) \cdot ( \zz -\bar{\zz}^0).
\end{align*}
For any $\zz \in B_{2\bar{r}} $,
\begin{equation}\label{est-vbar}
|\bar{v}(\zz)| \le  C  \bar{r}^{1+\ga}[ \nabla v_R ]_{0,\ga;B_{3\bar{r}} }.
\end{equation}
When $z_2^0 <1$ and $z_1^0 >k_1$, the Schauder interior estimate,
together with \eqref{est-vbar} and the $C^{1,\ga}$--estimate \eqref{est-DR1alpha},
leads to
\begin{equation}\label{est-vbar2ga}
\|\bar{v} \|'_{2,\ga; B_{\bar{r}}}
\le  C \|\bar{v} \|_{0,0;B_{2\bar{r}}}
\le CM (z_2^0)^{1+\ga}.
\end{equation}
Finally, for $\zz^0 \in \D^{k_1}$, by the corner estimate \eqref{est-corner1ga},
\[
[ \nabla v_R ]_{0,\ga;B_{3\bar{r}} } \le CM |\zz^0|^{-1}.
\]
Therefore, we have
\begin{equation}\label{est-vbar2gacorner}
\|\bar{v} \|'_{2,\ga; B_{\bar{r}}} \le  C \|\bar{v} \|_{0,0;B_{2\bar{r}}} \le CM |\zz^0|^{-1}(z_2^0)^{1+\ga}.
\end{equation}
Estimates \eqref{est-decay2ga2} and \eqref{est-vbar2ga}--\eqref{est-vbar2gacorner} imply
\begin{equation}
\|v_R\|_{2,\ga;(1+\gb,0);\D^{\frac{R}{4}}}^{(-\ga;\O)(-1-\ga;\mathcal{L}_1)} \le  CM.
\end{equation}
Taking $R=n$, we obtain a sequence $\{v_n\}_{n\in \N}$.
We can  choose a proper subsequence $\{v_{n_i}\}_{i \in \N}$
such that $\{v_{n_i}\}$ converges to $v$ in
$C_{(-\ga';\O)(-1-\ga';\mathcal{L}_1)}^{2,\ga';(1+\gb)}(\D^{\frac{n_i}{4}})$ for all $i\in \N$,
where $0< \ga' < \ga$.
Therefore, the limit function $v$ is a solution with estimate \eqref{est-v1}.

\subsection{Uniqueness of the solution}
Suppose that $v, \bar{v} \in  C_{(-\ga;  \O)(-1-\ga;\mathcal{L}_1)}^{2,\ga;(1+\gb,0)}(\D)$
both are the solutions for problem \eqref{eqn-v}--\eqref{con-g6}.
Then $\tv:= v-\bar{v}$ is also a solution of \eqref{eqn-v} with $g_5$ and $g_6$ vanishing in \eqref{con-g5}
and \eqref{con-g6}, respectively.
$\tv \in  C_{(-\ga;  \O)(-1-\ga;\mathcal{L}_1)}^{2,\ga;(1+\gb,0)}(\D)$ implies that
$|\tv(\zz)|$ decays as $|\zz| \to \infty$.
For any small $\ve>0$, there exists $R>0$ such that $|\tv(\zz)|< \ve$ on $L_R$.
Thus, by applying the maximum principle, we see that
$\|\tv\|_{0,0; \D^R} \le \ve$.
We know that $\tv\equiv 0 $ in $\D$ as $\ve \to 0$.

\subsection{Proof of part (ii) of Lemma \ref{thm-lemma} }
The procedure of proving part (ii) is primarily parallel to that of part (i),
except for the different regularity at the corner and decay rate
due to the opposite sign of $\nu_t$.
For the $C^{1,\ga}$--regularity at corner $\O$, we can estimate
$\bar{v} = v_R (\zz) - v_R(0,0) - \mathbf{a} \cdot \zz$,
where $\mathbf{a}=(a_1, a_2)$ is solved from the equations:
$$
\begin{cases}
a_1 = g_5'(0),\\[1mm]
\mathbf{a} \cdot  \boldsymbol{\nu} = g_6(0).
\end{cases}
$$
Since $\nu_n >0$ and $\nu_t > 0 $ imply that $\nu_2 >0$,
the equations above are uniquely solvable for $\mathbf{a}$.
Once we prove that $v_R$ is $C^{1,\ga}$ up to corner $\O$,
we can see that $\mathbf{a} = \nabla v_R(0,0)$.

We use
$$
v_3 = Mr^{1+\ga} \big(C_5 \sin((1+\ga +\tau)\t +\t_0) + C_6 (\sin \t)^{\ga}\big)
$$
to control $\bar{v}$ near corner $\O$.
In fact, denoting $\bar{\t}:=(1+\ga +\tau)\t +\t_0$, we have
\begin{align*}	
\Delta v_3 &= C_5M\big((1+\ga)^2 -(1+\ga +\tau)^2\big)r^{-1+\ga}\sin \bar{\t} \\
&\quad + C_6 M r^{-1+\ga} \big(\ga(-1+\ga)  (\sin \t)^{-2+\ga} + (1+2\ga) (\sin \t)^{\ga}\big)\\
&\le -C_5c_3 Mr^{-1+\ga} (\sin \t)^{-2+\ga}
\end{align*}
by choosing $\frac{C_5}{C_6}$ large enough.
Then we compute
\begin{align*}
&\big|(a_{ij} -\d_{ij})\partial_{ij} v_3\big|
 \le CMC_5\d r^{-1+\ga} (\sin \t)^{-2+\ga},\\[2mm]
&|(a_{ij})_{z_j} (v_3)_{z_i}|
 = O(\d)M (z_2+1)^{-1-\gb}r^{\ga} (\sin \theta)^{-1+\ga}
 \le C M\d r^{-1+\ga} (\sin \t)^{-2+\ga} .
\end{align*}
The estimates above yield
\begin{align*}
(a_{ij}(v_3)_{z_i})_{z_j}
= \big(\Delta +(a^\varphi_{ij} - \d_{ij})\partial^2_{z_i z_j}\big)v_3
+(a^\varphi_{ij})_{z_j} (v_3)_{z_i} \le 0
\end{align*}
for sufficiently small $\d$.

On the boundaries, we use that $\nu_n>0$, $\nu_t \ge 0$,
and $\bar\theta\in (0, \frac{\pi}{2})$ to obtain
\begin{align*}	
\left. \frac{Dv_3}{D \boldsymbol{\nu}}\right|_{\mathcal{L}_2}
&= Mr^{\ga} \Big\{C_5(\nu_n(1+\ga +\tau)\cos\bar{\t} + \nu_t (1+\ga) \sin \bar{\t}) \\
&\qquad\qquad  + C_6 \left. (\sin \t)^{-1+\ga}\big(\nu_t(1+\ga)\sin \t +  \nu_n\ga \cos \t\big)\Big\}\right|_{\t =\omega_0}\\[1mm]
& \ge  C_5 c_4  Mr^{\ga} \\
& \ge  g_6 -g_6(0)   \qquad\quad   (\text{by choosing suitably large $C_5$}),\\[2mm]
v_3|_{\mathcal{L}_1}& = C_5M z_1^{1+\ga} \sin \t_0
  \ge g_5(z_1) -g_5(0) -g_5'(0)z_1,\\[1mm]
v_3|_{L_R}& \ge  g_5(R) -g_5(0)-g_5'(0)R.
\end{align*}
Thus, by the comparison principle, we conclude
\begin{equation*}
|\bar{v}(\zz)| \le CM |\zz|^{1+\ga}.
\end{equation*}

On the other hand, the fact that $\nu_t\ge 0$ results in the decay rate $r^{-\gb}$,
which is slower than part (i) ($\nu_t < 0$). This can be achieved by setting
\[
v_4 = C M r^{-\gb}\big(\sin ( (\gb+\tau)\t + \t_0) + (\sin \t)^{\ga}\big).
\]
In the same way as in part (i), we can prove that $v_4$ is a supersolution of \eqref{eqn-v2ga}.
What is different from part (i) is that, for $\nu_t \ge 0$, we require $\gb$ small to guarantee
the positivity of $\frac{Dv_4}{D\boldsymbol{\nu}}$ on $\mathcal{L}_2$.
In fact, we have
\begin{align*}	
\left. \frac{Dv_4}{D \boldsymbol{\nu}}\right|_{\mathcal{L}_2}
= C Mr^{-1-\gb}
 &\big\{\nu_n\big((\gb +\tau)\cos((\gb +\tau)\omega_0 + \t_0)+ \ga (\sin \om_0)^{-1+\ga} \cos\om_0\big) \\
  &\,\,\, -\nu_t \gb \big(\sin ((\gb +\tau)\omega_0+ \t_0) + (\sin \om_0)^\ga\big)\big\},
\end{align*}
which is greater than $g_6$ if $\gb$ is small and $C$ is large.
After we obtain the $C^0$--estimate,
we apply the standard Schauder estimates
with proper scalings to achieve
estimate \eqref{est-v} in part (ii).

\section{Construction of the iteration map $\mathcal{Q}$}
\label{sec-iteration}
We first focus on \textbf{Problem WT}.

For a given upstream flow $U^-$ and $b$ in the slip condition \eqref{slipcon} satisfying
\[
\|U^- -U^-_0\|_{2,\ga;(1+\b,0);\D^-}  +\| b'\|^{(-\ga;0)}_{1,\ga; (1+\b);\R^+} \le \ve,
 \]
we define a map $\mathcal{Q}$ from $\Sigma^{C_0 \ve}$ to itself, provided that $C_0$ and $\ve$ are chosen properly,  where  $\Sigma^{C_0 \ve}$ is given as follows:
 \begin{align}\label{def-sapcesigma}
 \begin{split}
 \Sigma_1^\tau &:= \{ v: \| v\|_{2,\ga;(0,1+\gb);\D}^{(-\ga;\mathcal{L}_1)} + \|v_{z_1}\|_{1,\ga;(1+\gb,1);\D}^{(1-\ga;\mathcal{L}_1)}\le \tau \},\\[1mm]
  \Sigma_2^\tau &:= \{ v: \| v\|_{2,\ga;(1+\gb,0);\D}^{(-\ga;\O)(-1-\ga;\mathcal{L}_1)} \le \tau \},\\[1mm]
  \Sigma_3^\tau &:= \{ v: \| v\|_{2,\ga;(1+\gb);\R^+}^{(-\ga;0)} \le \tau \},\\[1mm]
  \Sigma^\tau  &:= \Sigma_1^\tau \times\Sigma_1^\tau\times\Sigma_2^\tau\times \Sigma_2^\tau\times\Sigma_3^\tau.
  \end{split}
 \end{align}
For notational convenience, we use $\|\cdot\|_{\Sigma_i}$ to denote the norm for $\Sigma_i^\tau$.
The norm  $\|\cdot\|_\Sigma$ is understood as the summation of the norms of all the components.
Given $V =(\d u_1, \d \rho, \d w , \d p, \d \hs') \in \Sigma^{C_0\ve}$,
we first solve equations \eqref{eqn-wp-lin1}--\eqref{eqn-wp-lin2}  with
the slip condition $\d\tilde{w}|_{\mathcal{L}_1} = b'$ and the boundary condition \eqref{con-g3}
on $\mathcal{L}_2$.
Once we obtain $(\d \tilde{w}, \d \tp)$, we use condition \eqref{con-g4} on $\mathcal{L}_2$
and equation \eqref{eqn-prho} to solve for $\d \tilde{\rho}$.
Then, by \eqref{eqn-u1}, we can compute $\d \tilde{u}_1$.
From equation \eqref{eqn-sigmaprime}, we update the shock function $\d \tilde{\s}$.
Thus, we can define $\mathcal{Q} (V) \equiv \widetilde{V}=(\d  \tilde{u}_1, \d \tilde{\rho}, \d \tilde{w},
\d \tp, \d \tilde{\s}' )$.

\subsection{Solve for $\d \tilde{w}$} \label{sec-solvedw}

We perform $\frac{\partial}{\partial z_2}$\eqref{eqn-wp-lin1}   $-\,  \frac{\partial}{\partial z_1}$\eqref{eqn-wp-lin2} to eliminate $\d \tp$ and obtain
 \begin{equation}\label{eqn-dw-z}
 \left( a_{ij}  \d \tilde{w}_{z_i}  \right)_{z_j}= 0,
\end{equation}
where
\begin{align*}
 a_{11} = \frac{(1- \d \hs'\gl_R)^2 +(\d \hs' \l_I)^2}{e \gl_I },\,\,\,
 a_{12} = a_{21}  = \frac{\gl_R - \d \hs'(\gl_R^2 + \l_I^2)}{e \gl_I},\,\,\,
 a_{22}= \frac{\gl_R^2 + \l_I^2}{e \gl_I } .
\end{align*}
In order to meet condition \eqref{con-aij} in Lemma \ref{thm-lemma},
we apply the following coordinate transformation:
\[
\begin{cases}
\bar{z}_1 = \sqrt{e^0 \l_I^0} z_1,\\[1mm]
\bar{z}_2= \sqrt{\frac{e^0}{\l^0_I }}z_2,
\end{cases}
\]
where $(e^0, \l_I^0)$ are $(e, \l_I)$ evaluated at the background state $U^+_0$.
Thus, equation \eqref{eqn-dw-z} becomes
 \begin{equation}\label{eqn-dw-zbar}
 \left( \bar{a}_{ij}  \d \tilde{\bar{w}}_{\bar{z}_i}  \right)_{\bar{z}_j}= 0,
\end{equation}
where
\begin{align*}
&\bar{a}_{11}(\bar{\zz}) =e^0 \l_I^0 a_{11}(\sqrt{\frac{1}{e^0 \l^0_I }}\bar{z}_1, \sqrt{\frac{\l^0_I}{ e^0}}\bar{z}_2),\quad
  \bar{a}_{22} (\bar{\zz})= \frac{e^0}{\l^0_I } a_{22}(\sqrt{\frac{1}{e^0 \l^0_I }}\bar{z}_1, \sqrt{\frac{\l^0_I}{ e^0}}\bar{z}_2),\\[1mm]
&\bar{a}_{12}(\bar{\zz}) = \bar{a}_{21} (\bar{\zz})= e^0 a_{12}(\sqrt{\frac{1}{e^0 \l^0_I }}\bar{z}_1, \sqrt{\frac{\l^0_I}{ e^0}}\bar{z}_2),
\quad
\d \tilde{\bar{w}}(\bar{\zz}) =\d \tilde{w}(\sqrt{\frac{1}{e^0 \l^0_I }}\bar{z}_1, \sqrt{\frac{\l^0_I}{ e^0}}\bar{z}_2).
\end{align*}
The boundary, $\mathcal{L}_2$, becomes $\bar{\mathcal{L}}_2: \bar{z}_1 = k_2 \bar{z}_2$ for $k_2 = k_1 \l^0_I$.
Condition \eqref{con-g3} becomes
\begin{equation}\label{con-g3bar}
\d \tilde{\bar{w}} + b_1 \d \tilde{\bar{p}} = \bar{g}_3
\end{equation}
in the $\bar{\zz}$--coordinates, where $\bar{g}_3$ is $g_3$ rescaled in the $\bar{\zz}$--coordinates.
Differentiating \eqref{con-g3bar} along  $\bar{\mathcal{L}}_2$ and using equations  \eqref{eqn-wp-lin1}--\eqref{eqn-wp-lin2}
to eliminate the $\d \tilde{\bar{p}}$ terms give rise to
\begin{equation}\label{con-g7}
\Big( k_2 -\frac{b_1}{e^0}(\bar{a}_{11} - k_2 \bar{a}_{12})\Big) (\d \tilde{\bar{w}})_{\bar{z}_1}
+ \Big(1+  \frac{b_1}{e^0}(k_2\bar{a}_{22} -  \bar{a}_{12}) \Big) (\d \tilde{\bar{w}})_{\bar{z}_2}  = \bar{g}'_3.
\end{equation}
Slightly modify \eqref{con-g7} into
\begin{equation}\label{con-g7bar}
\mu_1(\d \tilde{\bar{w}})_{\bar{z}_1} + \mu_2 (\d \tilde{\bar{w}})_{\bar{z}_2}  =  \bar{g}_7,
\end{equation}
where
\begin{align*}
\mu_1 & = k_2 -\frac{b_1}{e^0},\quad \mu_2 = 1+  \frac{b_1}{e^0}k_2,\\[1mm]
\bar{g}_7& =\bar{g}'_3+\frac{b_1}{e^0}\big((\bar{a}_{11}-1 )- k_2 \bar{a}_{12}\big) (\bar{w})_{\bar{z}_1}
 - \frac{b_1}{e^0}\big(k_2(\bar{a}_{22}-1)-\bar{a}_{12}\big)(\bar{w})_{\bar{z}_2}.
\end{align*}
Conditions \eqref{con-g7} and \eqref{con-g7bar} are equivalent when $\tilde{w}= w$,
{\it i.e.},  when  $V =(\d u_1, \d \rho, \d w , \d p, \d \hs')$ is a fixed point of $\mathcal{Q}$.
For \textbf{Problem WT}, $b_1 < 0$ (see Remark \ref{signb1}).
Then we normalize $\boldsymbol{\mu} = (\mu_1, \mu_2)$ into $\boldsymbol{\nu}=\frac{\boldsymbol{\mu}}{|\boldsymbol{\mu}|}$
and compute
\[
\nu_n = \frac{-b_1}{\sqrt{(e^0)^2+ b_1^2}} >0, \qquad \nu_t = \frac{-e^0}{\sqrt{(e^0)^2+ b_1^2}} <0.
\]
Moreover, we have
\[
\|\bar{a}_{ij}-\d_{ij}\|^{(-\ga;\mathcal{L}_1)}_{1,\ga;(0,1+\gb);\D} \le C C_0 \ve \le \d
\]
for sufficiently small $\ve$, so that condition \eqref{con-aij} is satisfied.
Therefore, applying  part  $\mathrm{ (i)} $ of  Lemma \ref{thm-lemma} and scaling back
to the $\zz$--coordinates, we have
\begin{eqnarray}
\| \d \tilde{w}\|^{(-\ga;  \O)(-1-\ga;\mathcal{L}_1)}_{2,\ga;(1+\gb,0);\D}
\le C \big(\|b'\|^{(-\ga;0)}_{1,\ga; (1+\b);\R^+}+\|g_7\|^{(1-\ga;0)}_{1,\ga; (2+\b);\R^+}\big),
\label{est-dw1}
\end{eqnarray}
where $g_7$ is $\bar{g}_7$ scaled back in the $\zz$--coordinates.
We know that
\begin{equation}
\| g_7 \|^{(1-\ga;0)}_{1,\ga; (2+\b);\R^+}
\le C\Big(\sum_{i=1,2}\| g_i\|^{(-\ga;0)}_{2,\ga; (1+\b);\R^+}+ \|V \|_{\Sigma}^2 \Big).
\label{est-g7}
\end{equation}
Since the Rankine-Hugoniot conditions \eqref{con-G1}--\eqref{con-G2} hold at the background states,
we have
\[
G_i(U_0^-, \bar{U}_0^+) =0, \qquad i=1,2.
\]
Therefore,  $g_i$ defined by \eqref{eqn-gi} can be rewritten as:
\begin{align*}
g_i = \grad_{\bar{U}} G_i(U_0^-, \bar{U}_0^+)\cdot \d \bar{U} - G_i(U_0^-, \bar{U})+ G_i(U_0^-, \bar{U}_0^+)
  +G_i(U_0^-, \bar{U}) - G_i(U_{\hs}^-, \bar{U}),
\end{align*}
which gives rise to the following estimates:
\begin{align}
\| g_i\|^{(-\ga;0)}_{2,\ga; (1+\b);\R^+}
&\le C \big(\|V\|_{\Sigma}^2+ \| \d U_{\hs}^- \|_{2,\ga;(1+\b,0);\D^-}\big) \nonumber \\
&\le C \big(\|V\|_{\Sigma}^2+ \| \d U^- \|_{2,\ga;(1+\b,0);\D^-} + \|\d \hs'\|_{\Sigma_3}\|\nabla  U^-\|_{1,\ga;(2+\b,0);\D^-}\big).\quad
\label{est-Usigma}
\end{align}
Combining  \eqref{est-g7} with \eqref{est-Usigma}, estimate \eqref{est-dw1} becomes
\[
 \| \d \tilde{w}\|_{\Sigma_2} \le C\big(1+ C_0^2 \ve  + C_0 \ve\big) \ve.
\]
Choosing $C_0 > 4C$ and $\ve < \frac{1}{C_0^2}$, we have
\begin{equation}\label{7.10a}
\| \d \tilde{w}\|_{\Sigma_2} \le 3C\ve < C_0 \ve,
\end{equation}
which implies that $\d\tilde{w}\in \Sigma^{C_0\ve}_2$.

\smallskip
\subsection{Higher decay rate for $(\d \tilde{w})_{z_1}$ } \label{set-extradecay}
In order to estimate the $C^0$--norm of $\d\tp$ in the next section, we need
an extra decay rate for $(\d \tilde{w})_{z_1}$ to control the logarithmic growth in $z_2$
({\it cf.} the argument from \eqref{eqn-dp} to \eqref{est-dpdecay}).

Differentiating \eqref{eqn-dw-z} with respect to $z_1$ yields
\begin{equation}\label{eqn-dwz1}
\big(a_{ij} (\d \tilde{w}_{z_1})_{z_i}\big)_{z_j}
=- \big((a_{ij})_{z_1} (\d\tilde{w})_{z_i}\big)_{z_j}.
\end{equation}
In domain $ \D^R\backslash D^{k_1}$, we solve the equation:
\begin{equation} \label{eqn-udw}
\big(a_{ij} u_{z_i}\big)_{z_j}
=- \big((a_{ij})_{z_1} (\d\tilde{w})_{z_i}\big)_{z_j} =: f,
\end{equation}
with the following Dirichlet boundary conditions:
\begin{align}
&u|_{\mathcal{L}_1 \cup  {\mathcal{L}_2} \cup  L_{k_1}} = (\d \tilde{w})_{z_1}|_{\mathcal{L}_1 \cup  {\mathcal{L}_2} \cup  L_{k_1}},\label{con-ul1l2}\\
&u|_{L_R}=  (\d \tilde{w})_{z_1} (R,0)
 + \big((\d \tilde{w})_{z_1} (R,\frac{R}{k_1})-(\d \tilde{w})_{z_1}(R,0)\big)\frac{k_1}{R} z_2.\label{con-ulr}
\end{align}
Condition \eqref{con-ulr} is artificially prescribed on $L_R$ so that the continuity of $u$ at the intersection points
of $L_R$ with $\mathcal{L}_1$ and ${\mathcal{L}_2}$ is achieved.

Given $R$, we obtain a solution $u_R$. The estimates of $u_R$ follow the same way as in Lemma \ref{thm-lemma}.
Once we have the desired {\it a priori} estimates, by the continuity method,
we also have the existence of the solution. Therefore, we only need
to point out the difference from the {\it a priori} estimates in Lemma \ref{thm-lemma}.

Equation \eqref{eqn-udw} with conditions \eqref{con-ul1l2}--\eqref{con-ulr} is a Dirichlet boundary problem with an inhomogeneous term
on the right-hand side. Notice that
\[
f  = O(\ve^2) |\zz|^{-2-2\gb}(z_2+1)^{-3} \big(\min(z_2,1)\big)^{-2+\ga},
\]
which implies
\begin{align}
|f| & \le C \ve^2 |\zz|^{-4-2\gb +\ga} (\sin \t)^{-2+\ga}\le \ve r^{-4-\gb } (\sin \t)^{-2+\ga},
\end{align}
provided that $\ga \le \gb$.

We use the barrier function
\[
v_5 = C \ve r^{-2-\gb} \big(\sin ((2+\gb+\tau)\t+ \t_0) + (\sin \t)^{\ga}\big),
\]
where $\gb, \t_0, \tau >0$ are small so that $(2+\gb+ \tau)\omega_0 + \t_0 <\pi$.
This can be achieved because $\omega_0 < \frac{\pi}{2}$.
It is easy to see that $u \le v_5$ on the boundary.
Following the same computation, we have
\begin{align*}
\big(a_{ij} (v_5)_{z_i}\big)_{z_j} \le -C\ve r^{-4-\gb} (\sin \t)^{-2\ga} \le f.
\end{align*}
Therefore, we conclude that
\begin{align*}
|u_R| \le C\ve  r^{-2-\gb}.
\end{align*}
With the $C^0$--estimate above, using the same scaling as in \S \ref{sec-keylemma},
we obtain the estimates in $\D^*:=\D^{\frac{R}{2}} \backslash \D^{k_1}$:
\begin{equation}
\|u\|_{1,\ga;(2+\gb,0);\D^*}^{(1-\ga;\O)(-\ga;\mathcal{L}_1)} \le C\ve.
\end{equation}
Choose a subsequence of $u_R$ so that, as $R\to \infty$, it converges to a solution $u$ of
\eqref{eqn-udw} in $\D\backslash \D^{k_1}$.
Since both $u$ and  $(\d \tilde{w})_{z_1}$ decay in the far field of domain $\D$,
the solution of  problem \eqref{eqn-udw}--\eqref{con-ulr} is unique.
Thus, we conclude that
\begin{equation}\label{decayu2gb}
| (\d \tilde{w})_{z_1}(\zz)| = |u(\zz)| \le C\ve |\zz|^{-2-\gb}.
\end{equation}

\smallskip
\subsection{Solve for $\d \tp$}

To solve for $\d \tp$,   we set the initial data for $\d \tp$  from  condition \eqref{con-g3}:
\begin{equation} \label{con-dp}
\d \tp = \frac{1}{b_1}(g_3 - \d \tilde{w})\qquad\,\, \mbox{on $\mathcal{L}_2$}.
\end{equation}
Using equation \eqref{eqn-wp-lin2}, we integrate in the $z_2$--direction to solve for $\d\tp$.
More precisely, let $\zz^0$ be any point in $\D $.
Let $\zz^I = (z_1^I, z_2^I)$ be the intersection point of $\mathcal{L}_2$ and the vertical line
passing through $\zz^0$.
By equation  \eqref{eqn-wp-lin2} with initial data  \eqref{con-dp},
we can express  $\d\tp$ explicitly in the following formula:
 \begin{equation} \label{eqn-dp}
 \d \tp (\zz^0) = \frac{1}{b_1} \big(g_3(z_2^I) - \d \tilde{w}(\zz^I)\big)
  + \int_{z_2^I}^{z_2^0} \big(-a_{11}(\d \tilde{w})_{z_1}-a_{12}(\d \tilde{w})_{z_2}\big)(z_1^0, s) ds.
 \end{equation}
 We first check the decay rate of $\d \tp$ by \eqref{eqn-dp}.
By the definition of $g_3$ in \eqref{5.14a} and estimate \eqref{est-Usigma} for $g_1$ and $g_2$,
together with estimate \eqref{7.10a} for $\delta \tilde{w}$, we have
 \[
 |g_3(z_2^I)-\d \tilde{w}(\zz^I)| \le C\ve |\zz^0|^{-1-\gb}.
 \]
 For the integral term in \eqref{eqn-dp},
 observe that $a_{12}=\frac{\gl_R - \d \hs'(\gl_R^2 +\l_I^2)}{e \gl_I }$, giving $ |\zz^0|^{-1-\gb}$ decay.
 Then we use \eqref{decayu2gb} and $\|\d \tilde{w}\|_{\Sigma_2} \le C\ve  $ to obtain
 \begin{align*}
 &\left|  \int_{z_2^I}^{z_2^0} \big(-a_{11}(\d \tilde{w})_{z_1}-a_{12}(\d \tilde{w})_{z_2}\big)(z_1^0, s) ds\right|\\
 &\le {} C\ve  |\zz^0|^{-1-\gb} \int^{z_2^I}_{z_2^0}|\zz^0|^{-1} \ ds
 \le {}  C\ve |\zz^0|^{-1-\gb} \frac{k_1 z_1^0 }{ |\zz^0|}\\
 &\le {}  C\ve |\zz^0|^{-1-\gb} .
 \end{align*}
 Therefore, we have
 \begin{equation} \label{est-dpdecay}
 | \d \tp (\zz^0)  |  \le C \ve |\zz^0|^{-1-\gb}.
 \end{equation}

\smallskip
 For the corner regularity, for any $\zz^0 \in \D^{k_1}$, equation \eqref{eqn-dp} implies
 \begin{equation}  \label{est-dp1alpha}
 | \d \tp (\zz^0) -  \d \tp (0,0) |  \le C \ve |\zz^0|^{ \ga},
 \end{equation}
 indicating that $ \d \tp $ is $C^{\ga}$ smooth up to corner $\O$.
 The estimates for the derivatives of $\d \tp$ follow from the observation below.

Recall that \eqref{eqn-dw-z} is obtained by
differentiation $\frac{\partial}{\partial z_2}$\eqref{eqn-wp-lin1}$ -\, \frac{\partial}{\partial z_1}$\eqref{eqn-wp-lin2}.
Notice that $\delta \tilde{w}$ satisfies \eqref{eqn-dw-z} and $(\delta \tilde{w}, \delta \tilde{p})$
satisfies \eqref{eqn-wp-lin2} in domain $\mathbb{D}$,
since $\delta \tilde{p}$ is solved from \eqref{eqn-wp-lin2} with initial
data \eqref{con-dp} (see \eqref{eqn-dp} for the expression
for $\delta \tilde{p}$).
Therefore, we obtain $\frac{\partial}{\partial z_2}\eqref{eqn-wp-lin1}$, {\it i.e.},
\begin{equation}\label{new-N}
\frac{\partial}{\partial z_2}((\delta \tilde{p} )_{z_1} -
a_{12}(\delta \tilde{w})_{z_1} - a_{22}(\delta \tilde{w})_{z_2}) =0.
\end{equation}
To recover \eqref{eqn-wp-lin1}, we integrate equation \eqref{new-N} along the $z_2$--direction to deduce
\begin{equation}\label{diff}
 (\d  \tilde{p})_{z_1} - a_{12}(\d \tilde{w} )_{z_1}  - a_{22} (\d \tilde{w} )_{z_2} = f(z_1),
 \end{equation}
where
$$
f(z_1) = ((\delta \tilde{p})_{z_1}- a_{12}(\delta \tilde{w})_{z_1}
  - a_{22}(\delta \tilde{w})_{z_2}) (z_1, \frac{z_1}{k_1}).
$$
Notice that condition \eqref{con-g7bar} is a modification from \eqref{con-g7}, so that $f$ does not vanish on $\mathcal{L}_2$.
 Using conditions \eqref{con-g3} and \eqref{con-g7bar}, together with the fact that equation \eqref{eqn-wp-lin2} holds
 up to boundary $\mathcal{L}_2$, we obtain
 \begin{align}\label{exp-f}
 f(z_1)= {}&\Big(\big(\frac{a_{11}}{k_1}- \frac{1}{k_1e^0 \l_I^0}- a_{12}\big)(\tilde{w}-w)_{z_1}
 -\big(a_{22} - \frac{\l_I^0}{e^0} -\frac{a_{12}}{k_1} \big)(\tilde{w}-w)_{z_2} \Big) (z_1, \frac{z_1}{k_1}).
 \end{align}

Equation \eqref{eqn-wp-lin1} will be recovered later, when we obtain a fixed point for $\mathcal{Q}$.
For now, we can use equations  \eqref{eqn-wp-lin2} and \eqref{diff} to estimate the derivatives
of $\d \tp$ in terms of $\d\tilde{w}$.
Thus, together with estimates \eqref{est-dpdecay}--\eqref{est-dp1alpha},
we see that $\d \tp \in \Sigma^{C_0\ve}_2$, by choosing large enough $C_0$.

\subsection{Solve for $(\d \tilde{\rho}, \d  \tilde{u}_1)$}
We use \eqref{con-g4} as the initial data on $\mathcal{L}_2$
and solve equation \eqref{eqn-prho} to obtain $\d \tilde{\rho}$ and directly compute $\d\tilde{u}_1$
by \eqref{eqn-u1}.
Since  $(\d\tilde{\rho}, \d\tilde{u}_1)$ are obtained by the algebraic equations,
it is obvious that the smoothness of  $(\d\tilde{\rho}, \d\tilde{u}_1)$ is the same as
that of $(\d \tilde{w}, \d \tp)$.
However, in equations \eqref{eqn-u1}--\eqref{eqn-prho}, both $\frac{p}{\rho^\g}$ and $B$ are conserved,
rendering the non-decay  of  $(\d \tilde{\rho}, \d  \tilde{u}_1)$ in the $z_1$--direction.
On the other hand,  $(\d \tilde{\rho}, \d  \tilde{u}_1)$ have the same decay rate as their initial data
on $\mathcal{L}_2$ in the $z_2$--direction.

More precisely, for any point $\zz \in \D$,
let $\zz^I $ be the intersection of $\mathcal{L}_2$ and the horizontal line passing through $\zz$.
Since $\frac{p}{\rho^{\gamma}}$ is constant along the $z_2$--direction,
we use
$$
\frac{\tilde{p}}{\tilde{\rho}^{\gamma}}(\mathbf{z}) =
 \frac{\tilde{p}}{\tilde{\rho}^{\gamma}}(\mathbf{z}^I)
$$
 to solve for $\delta \tilde{\rho}$:
 \begin{align*}
\d \tilde{\rho}(\zz)
={}&\Big( \frac{\tp(\zz)}{\tp(\zz^I)} \Big)^{\frac{1}{\g}} \rho (\zz^I) -{\rho}^+_0  \\
={}&\Big( \frac{\tp(\zz)}{\tp(\zz^I)} \Big)^{\frac{1}{\g}}
 \big(g_4(z_2) -b_2 \d \tilde{w}(\zz^I) -b_3 \d \tp (\zz^I)\big)
 +\Big( \big( \frac{\tp(\zz)}{\tp(\zz^I)} \big)^{\frac{1}{\g}} - 1 \Big) \rho^+_0,
\end{align*}
where $\tilde{p} = p_0^+ +\delta \tilde{p}$ and $\tilde{\rho} = \rho_0^+ + \delta \tilde{\rho}$.
From the above expression, we see that $|\d \tilde{\rho} (\zz)| \le C \ve (z_2 +1) ^{-1-\gb}$.
The derivatives of $\d \tilde{\rho} $ also decay with appropriate rate adapted to the corresponding norms
in the  $z_2$--direction.
Thus, we have
\[
\|\d \tilde{\rho}\|_{2,\ga;(0,1+\gb);\D}^{(-\ga;\mathcal{L}_1)} \le C\ve.
\]

To see that $\d \tilde{\rho} \in \Sigma_1^{C_0\ve}$,
we need to obtain the estimate for the other part in the norm ({\it cf.} \eqref{def-sapcesigma}).
For this purpose, we rewrite the expression of  $\d \tilde{\rho}$ into the following form:
\begin{equation}\label{exp-drho}
\d \tilde{\rho}(\zz) = A(z_2)\tp(\zz)^{\frac{1}{\g}}- \rho^+_0.
\end{equation}
Taking the partial derivative with respect to $z_1$ on \eqref{exp-drho} yields
\[
(\d \tilde{\rho})_{z_1} = \frac{1}{\g}A(z_2)\tp(\zz)^{\frac{1}{\g}-1}(\d \tp)_{z_1}.
\]
The expression above shows that $(\d \tilde{\rho})_{z_1}$ and $(\d \tp)_{z_1}$ have the same decay pattern,
giving the estimate:
\[
\|(\d \tilde{\rho})_{z_1}\|_{1,\ga;(1+\gb,1);\D}^{(1-\ga;\mathcal{L}_1)}\le C\ve
.\]
The same argument also applies to the decay of $\d\tilde{u}_1$.
Thus, we conclude that $\d\tilde{\rho}, \d\tilde{u}_1 \in \Sigma^{C_0\ve}_1$.

\subsection{Update shock-front}  \label{sec-updateshock}

From \eqref{eqn-sigmaprime},  we can update $\d \tilde{\s}' $ by
\begin{equation} \label{eqn-sigmapt}
\tilde{ \s}' (z_2) = \frac{[ \tilde{u}_1 \tilde{w} ]}{[\tp]}(k_1 z_2, z_2),
\end{equation}
where the left state is $U^-_{\hs}$.

To estimate  $\d \tilde{ \s}' $, first let
 \[G(U^- ,U) := \frac{u_1 w - u_1^- w^-}{p -p^-}.
 \]
 Then equation \eqref{eqn-sigmapt} can be written as
 \begin{equation}\label{eqn-sigma2}
\tilde{ \s}' (z_2) =G(U_{\hs}^- ,\widetilde{U})(k_1 z_2, z_2).
\end{equation}
We know that \eqref{eqn-sigma2} is satisfied for the background states so that
\begin{equation}\label{eqn-sigma0}
\hs_0' (z_2) = k_1=G(U_0^- ,U_0^+).
\end{equation}
Taking the difference between equations \eqref{eqn-sigma2} and  \eqref{eqn-sigma0} gives
\begin{equation}\label{eqn-dsigma}
\d  \tilde{ \s}' (z_2) =G(U_{\hs}^- ,\widetilde{U})(k_1 z_2, z_2) -G(U_0^- ,U_0^+),
\end{equation}
which gives rise to the following estimates, similar to \eqref{est-Usigma}:
\begin{align*}
\|\d \tilde{\s}'\|_{\Sigma_3}
& \le  C \big(\| \d U_{\hs}^- \|_{2,\ga;(1+\b,0);\D^-}+\|\d \tilde{u}_1\|_{\Sigma_1} + \|(\d \tilde{w}, \d \tp)\|_{\Sigma_2}\big)\\
 &\le C\big(\| \d U^- \|_{2,\ga;(1+\b,0);\D^-} +  \| \d \hs' \|_{\Sigma_3}\|\nabla  U^- \|_{1,\ga;(2+\b,0);\D^-} + \ve\big)\\
&\le C\ve.
\end{align*}
Choosing  $C_0 > C$, we see that $\d  \tilde{ \s}'  \in \Sigma_3^{C_0\ve} $.
Therefore, we construct a map $\mathcal{Q}$ from $\Sigma^{C_0\ve}$ to itself.

\subsection{Fixed point of $\mathcal{Q}$}
 We use the Schauder fixed point theorem to prove the existence of the subsonic solution and the transonic shock.
 To fit into the framework of the Schauder fixed point theorem, we define the following Banach space:
 \[
 \Sigma':= \big\{(f_1,f_2,f_3,f_4,f_5): \|(f_1,f_2)\|_{\Sigma'_1}+ \|(f_3,f_4)\|_{\Sigma'_2}+\|f_5\|_{\Sigma'_3} < \infty\big\},
 \]
 where $\|\cdot \|_{\Sigma'_i}, i=1,2,3$, are the same norms defined in \eqref{def-sapcesigma},
 except that $\ga$ is replaced by $\ga'$, where $0 <\ga' <\ga$.
 Thus, $\Sigma^{C_0\ve}$ is a nonempty, convex, and compact subset of $\Sigma'$,
 and $\mathcal{Q}$ is a map from $\Sigma^{C_0\ve}$ into itself.
 Once we can show that $\mathcal{Q}$ is continuous, by the Schauder fixed point theorem, there is  a fixed point of $\mathcal{Q}$.
 To show the continuity of $\mathcal{Q}$, we can use the following compactness argument.

 On the contrary, assume that $\mathcal{Q}$ is not continuous. Then there exist a sequence $\{V^n\}_{n\in \N}$, a function $V^0$ in $\Sigma^{C_0\ve}$,
 and a constant $\d_0 >0$ such that $V^n \to V^0$ in $\Sigma'$,
 while $\|\mathcal{Q}V^n -\mathcal{Q}V^0\|_{\Sigma'} \ge \d_0$.
 Since $\{\mathcal{Q}V^n\}\subset \Sigma^{C_0\ve}$, which is compact in $\Sigma'$,
 we can select a subsequence $\mathcal{Q}V^{n_k}$ such that $\mathcal{Q}V^{n_k}\to W^0 \in \Sigma^{C_0\ve}$ as $k\to \infty$.
 Following the iteration process in \S \ref{sec-solvedw}--\S \ref{sec-updateshock},
 we see that $W^0 = \mathcal{Q}V^0$, which leads to a contradiction.
 This shows that $\mathcal{Q}$ is a continuous map from $\Sigma^{C_0\ve}$ into itself.

 Therefore, by the Schauder fixed point theorem, there exists a fixed point of $\mathcal{Q}$,
 denoted by $V=(\d u_1, \d \rho, \d w, \d p, \d \hs')$.
 Thus, $U =(u_{10}^+ +\d u_1, (u_{10}^+ +\d u_1)\d w, \rho_0^+ +\d \rho, p_0^+ + \d p)$ gives a subsonic solution,
 and $\hs$ gives the transonic shock-front in the $\yy$--coordinates.
 Therefore, we have proved the existence of solutions in  part (i) of Theorem \ref{thm-lag}.

\section{Uniqueness of the transonic solutions}  \label{sec-unique}

Let $V^i = (\d u^i_1, \d \rho^i, \d w^i , \d p^i, {\d \hs^i}^{\prime} ) \in \Sigma^{C_0\ve}, i=1,2$, be two fixed points of $\mathcal{Q}$.
Set
$$
V^d =(\d u^d_1, \d \rho^d, \d w^d, \d p^d, {\d \hs^d}^{\prime} ) = V^2- V^1.
$$

Denote $w^i$ scaled in the $\bar{\zz}$--coordinates by $\bar{w}^i$, and the rest of the variables are denoted in the same manner.
By the construction of $\mathcal{Q}$, we know that $\d\bar{w}^i$ satisfies \eqref{eqn-dw-zbar} for $i=1,2$.
Then taking the difference of the two equations results in
\begin{equation} \label{eqn-wd}
\big(\bar{a}_{ij}(\bar{V}^2) (\d {\bar{w}}^d)_{\bar{z}_i}\big)_{\bar{z}_j}
= -\big((\bar{a}_{ij}(\bar{V}^2)-\bar{a}_{ij}(\bar{V}^1))(\d w^1)_{\bar{z}_i}\big)_{\bar{z}_j} =:\bar{f}.
\end{equation}

The inhomogeneous term $\bar{f}$ in \eqref{eqn-wd} will result in the lower decay rate for $\d{\bar{w}}^d$.
In definition \eqref{def-sapcesigma}, we replace $\gb$ with $\frac{\gb}{2}$ and denote the new norms
by $\|\cdot\| _{\tgs_i}, i=1,2,3$.

Set $M_1 = \|V^d \|_{\tgs}$.
If $M_1 = 0$, we see that $V^1 =V^2$.
Now suppose that $M_1 >0$.
Then the estimates follow the same way as in the proof of Lemma \ref{thm-lemma},
except that we need to take care of the inhomogeneous term $\bar{f}$.

We first estimate $\bar{f}$ as follows: For $|\bar{\zz}|\ge 1$,
\begin{align*}
|\bar{f}(\bar{\zz})|& \le C\ve \|V^d \|_{\tgs} |\bar{\zz}|^{-1-\gb} \big(\max(\bar{z}_2,1)\big)^{-2-\frac{\gb}{2}}\big(\min(\bar{z}_2,1)\big)^{-1+\ga}\\
& \le C M_1 \ve  |\bar{\zz}|^{-1-\gb} \bar{z}_2^{-2+\ga}\\
& =  C M_1 \ve  |\bar{\zz}|^{-3+\ga-\gb} (\sin \t)^{-2+\ga}\\
& \le  C M_1 \ve  |\bar{\zz}|^{-3-\frac{\gb}{2}} (\sin \t)^{-2+\ga},
\end{align*}
provided that $\ga \le \frac{\gb}{2}$.
For  $|\bar{\zz}| < 1$,
\begin{align*}
|\bar{f}(\bar{\zz})|& \le C C_0\ve \|V^d \|_{\tgs} |\bar{\zz}|^{ -1} \bar{z}_2^{-1+\ga}
\le  CC_0\ve M_1 r^{-2+\ga} (\sin \t)^{-1+\ga}.
\end{align*}

Then we use the barrier function $v_6 $, similar to $v_1$ in \S \ref{C0estimates} for the corner estimates:
\[
v_6= C_7 M_1 \ve \big(r^{\ga} \sin(( \ga + \tau)\t + \t_0  ) +\bar{z}_2^{\ga}\big).
\]
Observe that
\begin{align*}
\big(\bar{a}_{ij}(\bar{V}^2)(v_6)_{\bar{z}_i}\big)_{\bar{z}_j}
\le -C_7 c_5M_1 \ve r^{-2+\ga} (\sin \t)^{-2+\ga} \le \bar{f},
\end{align*}
when $C_7$ is chosen large enough.

Since $\d\bar{w}^d$ vanishes on $\mathcal{L}_1$, then we use \eqref{con-g7} to obtain
\begin{align*}
\left.\frac{ D( \d \bar{w}^d)}{D \boldsymbol{\nu}}\right|_{\bar{\mathcal{L}}_2}={}& b_4 g^d_7(\bar{z}_2)\\
 = {}& b_4\frac{d}{d\bar{z}_2}\left(\left(\bar{g}_3(U^-_{\hs^2}, \bar{U}^2)-\bar{g}_3(U^-_{\hs^1}, \bar{U}^1)\right)(k_2\bar{z}_2,\bar{z}_2)\right) \\[1mm]
 & + \bar{g}_8 (\bar{V}^2, \grad_{\bar{\zz}} \bar{w}^2)- \bar{g}_8 (\bar{V}^1, \grad_{\bar{\zz}} \bar{w}^1),
\end{align*}
where
\begin{align*}
 &b_4 ={}  \frac{e^0}{\sqrt{((e^0)^2 + b_1^2)(k_2^2 +1 )}},\\[1mm]
&\bar{g}_8(\bar{V}, \grad_{\bar{\zz}} \bar{w})  ={}\frac{b_1}{e^0}\big((\bar{a}_{11}(\bar{V}) -1)- k_2 \bar{a}_{12}(\bar{V})\big)\bar{w}_{\bar{z}_1}
- \frac{b_1}{e^0}\big(k_2(\bar{a}_{22}(\bar{V})-1) -  \bar{a}_{12}(\bar{V})\big)\bar{w}_{\bar{z}_2}.
\end{align*}

Thus, we have the following estimates for $g_7^d$:
\begin{align*}
 \|g^d_7\|_{1,\ga;(2+\frac{\gb}{2});\R^+}^{(1-\ga;0)}
 &\le{}\sum_{i=1,2}C\Big(\big\|\big(g_i(U^-_{\hs^2}, \bar{U}^2)-g_i(U^-_{\hs^1},\bar{U}^1)\big)\big|_{\bar{\mathcal{L}}_2}\big\|^{(-\ga;0)}_{2,\ga; (1+\frac{\b}{2});\R^+}
    + C_0\ve \|V^d\|_{\tgs}  \Big)\\
&\le{} C\big( C_0\ve \|V^d\|_{\tgs}+ \| \d U_{\hs^2}^-  -  \d U_{\hs^1}^- \|_{2,\ga;(1+\frac{\gb}{2},0);\D^-}\big)\\
&\le{} C\big(C_0\ve  \|V^d\|_{\tgs} + \| \d U^- \|_{2,\ga;(1+\frac{\gb}{2},0);\D^-} \|{\d \hs^d}^{\prime}\|_{\tgs_3}\big) \\
&\le{} C C_0\ve M_1,\\[1mm]
g^d_7 (\bar{z}_2) \le {} &CC_0\ve M_1 \big(\max(\bar{z}_2,1)\big)^{-2-\frac{\gb}{2}}\big(\min(\bar{z}_2,1)\big)^{-1+\ga} \\
\le {}  &CC_0\ve M_1\big(\max(r,1)\big)^{-2-\frac{\gb}{2}} \big(\min(r,1)\big)^{-1+\ga}.
\end{align*}
Thus, we conclude
\begin{align*}	
\left. \frac{Dv_6}{D \boldsymbol{ \nu}} \right|_{\bar{\mathcal{L}}_2} \ge {}& C  c_2M_1 \ve r^{-1+\ga}
\ge  b_4 g^d_7(\bar{z}_2) = \left.\frac{ D( \d \bar{w}^d)}{D \boldsymbol{\nu}}\right|_{\bar{\mathcal{L}}_2}.
\end{align*}

On the cutoff boundary $L_R$, we know
\begin{align*}
\d \bar{w}^d & \le C C_0 \ve R^{-1-\gb} \le  CM_1\ve R^\ga \le v_6
\end{align*}
for sufficiently large $R$.
Therefore, we can use $v_6$ to bound $\d \bar{w}^d$ in $\D^R$. For the decay in $\D^R\backslash \D^{k_1}$, we use
\[
v_7 (\bar{\zz})= M_1 \ve r^{-1-\frac{\gb}{2}}\big(C_3 \sin ((1+\tfrac{\gb}{2} +\tau)\t + \t_0)+C_4 (\sin \t)^{\ga} \big).
\]

The same calculation as in \S \ref{C0estimates} shows that
\begin{align*}
\big(\bar{a}_{ij}(\bar{V}^2)(v_7)_{\bar{z}_i}\big)_{\bar{z}_j}& \le -C M_1 \ve r^{-3-\frac{\gb}{2}} (\sin \t)^{-2+\ga} \le \bar{f}.
\end{align*}
It is also easy to verify that
\begin{align*}	
\left. \frac{Dv_7}{D \boldsymbol{ \nu}} \right|_{\bar{\mathcal{L}}_2}
&\ge  CM_1\ve r^{-2-\frac{\gb}{2}}  \ge \left.\frac{ D( \d \bar{w}^d)}{D \boldsymbol{\nu}}\right|_{\bar{\mathcal{L}}_2}.
\end{align*}
We choose $R$ large enough, so that $R^{-\frac{\gb}{2}} \le  M_1$. Therefore, we obtain the control on  $L_R$:
\begin{align*}
\d \bar{w}^d & \le C C_0 \ve R^{-1-\gb} \le  C C_0 M_1\ve R^{-1-\frac{\gb}{2}} \le v_7.
\end{align*}
By the comparison principle, we conclude
\begin{equation*}
|\d \bar{w}^d (\bar{\zz})| \le CM_1 \ve |\bar{\zz}|^{-1-\frac{\gb}{2}}   \qquad  \mbox{ for }  \zz \in \D^R \backslash \D^{k_1}.
\end{equation*}
Once we have the $C^0$--estimates above, the rest is similar to the procedure as in \S \ref{sec-keylemma}.
In the end, we have
\[
\|V^d\|_{\tgs} \le C\ve M_1 = C\ve \|V^d\|_{\tgs}.
\]
Choose $\ve$ sufficiently small, so that $ C\ve <\frac{1}{2}$.
We see that $M_1=0$, which contradicts our assumption that $M_1>0$.
This completes the proof of the uniqueness of the solution
for \textbf{Problem WT} in Theorem \ref{thm-lag}.

\section{Asymptotic behavior of the subsonic solution}  \label{sec-asymptotic}

The estimate that  $\|V\|_{\Sigma} \le C_0 \ve$ implies
\[
\|\d p\|_{2,\ga; (1+\gb,0);\D}^{(-\ga;\O)(-1-\ga;\mathcal{L}_1)} \le C_0 \ve,
\qquad \|\d w\|_{2,\ga; (1+\gb,0);\D}^{(-\ga;\O)(-1-\ga;\mathcal{L}_1)} \le C_0 \ve.
\]
This means that $p \to p_0^+$  and $\frac{u_2}{u_1} \to 0$ at rate $|\zz|^{-1-\gb}$.
However, for fixed $z_2$, $(u_1, \rho)$ does not converge to $(u_{10}^+, \rho_0^+)$,
as $z_1 \to \infty$.
Observe that, from \eqref{exp-drho}, $\rho$ can be expressed by
\begin{equation}\label{exp-rho}
\rho (\zz) = A(z_2) p(\zz)^{\frac{1}{\g}},
\end{equation}
where $A$ can be solved from the Rankine-Hugoniot conditions \eqref{con-RH1}--\eqref{con-RH4}
when we find the shock function $\hs$. Then we define the limit for $\rho$ in the far field:
\begin{equation}\label{exp-rhoinfty}
\rho^\infty (z_2) = A(z_2) (p^+_0)^{\frac{1}{\g}}.
\end{equation}
 Taking the difference between \eqref{exp-rho} and \eqref{exp-rhoinfty} yields
 \begin{align*}
 \|\rho -\rho^\infty\|_{2,\ga;(1+\gb,0);\D}^{(-\ga;\mathcal{L}_1)}
 \le C \| \d p\|_{2,\ga;(1+\gb,0);\D}^{(-\ga;\O)(-1-\ga;\mathcal{L}_1)} \le C_0\ve.
 \end{align*}
In the same way, we use \eqref{eqn-u1} to obtain the limit for $u_1$:
\[
u_1^\infty(z_2) = \sqrt{2B(z_2) - \frac{2\g p_0^+ }{ (\g -1 )\rho^\infty(z_2) }}.
\]
Similarly, we have
\begin{align*}
 \|u_1 -u_1^\infty\|_{2,\ga;(1+\gb,0);\D}^{(-\ga;\mathcal{L}_1)} & \le C_0\ve.
\end{align*}

Since $\d{\hs}^{\prime}\in\Sigma_3^{C_0\ve}$,  the coordinate transformation \eqref{eqn-tranz}
has higher regularity than $U$ in the $\zz$--coordinates.
Therefore, $\|V\|_{\Sigma}\le C_0 \ve$  with the estimates above yields the corresponding
estimate \eqref{est-U-small-pert3} in the $\yy$--coordinates.
Thus, we have proved part (i) of  Theorem \ref{thm-lag}.

The coordinate transformation \eqref{def-coord} between $\xx$ and $\yy$ also has
higher regularity than $U$ and $\hs'$ in the subsonic domain $\D^{\hs}$,
and is bi-Lipschitz across the shock-front $\mathcal{T}$,
thanks to the Rankine-Hugoniot conditions \eqref{con-RH1}.
Therefore, estimate \eqref{est-U-small-pert3} implies estimate \eqref{est-U-small-pert},
so that the proof of part (i) of Theorem \ref{thm-main} is completed.

\section{Key points in solving  \textbf{Problem ST}} \label{sec-th}
Based on part (ii) of Lemma \ref{thm-lemma}, we can prove
part (ii) of  Theorem \ref{thm-lag}  in the same way as above.
Since most of the proof is parallel to that in part (i) of  Theorem \ref{thm-lag},
we will only point out the difference from part (i).

\subsection{Estimates for the existence of  solutions}
The procedure to construct the iteration map is the same as in part (i)
for \textbf{Problem WT}.
In \S \ref{set-extradecay},  we obtain the faster decay rate, $r^{-2-\gb}$, for $(\d\tilde{w})_{z_1}$,
compared to the $r^{-1-\gb}$ decay for $\d \tilde{w}$.
For  \textbf{Problem ST}, we can only gain extra $\frac{\gb}{2}$ decay rate, {\it i.e.},
$r^{-\frac{3\gb}{2}}$ decay for  $(\d \tilde{w})_{z_1}$.
Specifically, we use
\[
v_8 = C M r^{-\frac{3\gb}{2}} \big(\sin((\tfrac{3\gb}{2}+\tau)\t + \t_0) + (\sin \t)^{\ga}\big)
\]
as the barrier function for  $(\d \tilde{w})_{z_1}$ to obtain the following estimate:
\begin{align*}
\big(a_{ij} (v_8)_{z_i}\big)_{z_j} \le -C\ve r^{-2-\frac{3\gb}{2}} (\sin \t)^{-2+\ga}.
\end{align*}
On the other hand, $f$ in \eqref{eqn-udw} satisfies
\begin{align*}
|f (\zz)| &= O(\ve^2) |\zz|^{-2\gb}\big(\max(z_2,1)\big)^{-3} \big(\min(z_2,1)\big)^{-1+\ga}\\[1mm]
 & \le C \ve^2 r^{-2-2\gb +\ga} (\sin \t)^{-2+\ga} \\[1mm]
& \le \ve r^{-2-\frac{3\gb}{2} } (\sin \t)^{-2+\ga},
\end{align*}
provided that $\ga \le \frac{\gb}{2} $, which is the same
restriction on $\ga$ and $\gb$ as in \S \ref{sec-unique}.

There is no difference for the boundary estimates. Thus, we conclude
\[
|(\d \tilde{w})_{z_1}| \le C\ve r^{-\frac{3\gb}{2}}.
\]
With the estimate above and using expression \eqref{eqn-dp} for $\d\tp$,  we see that
\begin{align*}
| \d \tp (\zz^0)  | \le{}& |g_3(z_2^I) - \d \tilde{w}(\zz^I)|
+  \left|  \int_{z_2^I}^{z_2^0} \big(-a_{11}(\d \tilde{w})_{z_1}-a_{12}(\d \tilde{w})_{z_2}\big)(z_1^0, s) ds   \right|\\
\le {}& C\ve  |\zz^0|^{-\gb} +  C\ve  |\zz^0|^{-\gb}\int^{z_2^I}_{z_2^0}|\zz^0|^{-\frac{\gb}{2} } (1+ s)^{-1} \ ds\\
\le {} &  C\ve |\zz^0|^{-\gb} +  C\ve  |\zz^0|^{-\gb}\int^{z_2^I}_{0} (1+ s)^{-1-\frac{\gb}{2}} \ ds\\
\le {} &  C\ve |\zz^0|^{-\gb} .
\end{align*}

\subsection{Uniqueness of subsonic solutions}

We need to take care of the decay estimate, since the rest is similar to those in \S \ref{sec-unique}.

Now $\bar{f}$ in \eqref{eqn-wd} can be controlled as follows:
\begin{align*}
|\bar{f}(\bar{\zz})|& \le C\ve \|V^d \|_{\tgs} |\bar{\zz}|^{-\gb} \big(\max(\bar{z}_2,1)\big)^{-2-\frac{\gb}{2}}\\[1mm]
& \le C M_1 \ve  |\bar{\zz}|^{-\gb} (\bar{z}_2)^{-2+\frac{\gb}{2}}\\[1mm]
& \le  C M_1 \ve  |\bar{\zz}|^{-2-\frac{\gb}{2}} (\sin \t)^{-2+\frac{\gb}{2}}.
\end{align*}

The barrier function
\[
v_9 = CM_1 \ve r^{-\frac{\gb}{2} } \big(\sin ( (\tfrac{\gb}{2}+\tau)\t + \t_0) + (\sin \t)^{\frac{\gb}{2}}\big)
\]
can be estimated as
\begin{align*}
\big(\bar{a}_{ij}(\bar{V}^2)(v_9)_{\bar{z}_i}\big)_{\bar{z}_j}
\le -C M_1 \ve r^{-2-\frac{\gb}{2}} (\sin \t)^{-2+\frac{\gb}{2}} \le \bar{f}.
\end{align*}
With similar boundary estimates, we can conclude the uniqueness of the subsonic solution.

Therefore, we have proved that, given a constant transonic flow on arc $\wideparen{TH}$ or $ \wideparen{TS}$,
if the upstream flow and the wedge boundary are perturbed,
then there exist a unique subsonic solution and transonic shock,
which are close to the background constant state and straight shock front.
This shows the stability of the constant transonic flows past wedges.
For the constant states on $ \wideparen{TS}$,
the regularity of the subsonic solution near the corner is $C^\ga$ and
the decay rate in the far field is $r^{-1-\gb}$.
Furthermore, we gain the higher decay rate $r^{-2-\gb}$ for the directional derivative
along the streamlines of $w =\frac{u_2}{u_1}$,
the direction of the flow.
On $\wideparen{TH}$, we obtain the $C^{1,\ga}$--regularity at the corner and $r^{-\gb}$ decay
in the far field.

\subsection*{Acknowledgments}
The research of Gui-Qiang G. Chen was supported in
part by
the UK
Engineering and Physical Sciences Research Council Award
EP/L015811/1
and the Royal Society--Wolfson Research Merit Award (UK).
The research of Jun Chen was supported in part by  Sustech Start-up Grant	Y01286205.
The research of Mikhail Feldman was
supported in part by the National Science Foundation under Grant
DMS-1401490 and the Van Vleck Professorship Research Award by the University
of Wisconsin-Madison.

\end{document}